\documentclass[11pt]{article}
\usepackage{amsmath,amssymb,euscript}

\makeatletter \@addtoreset{equation}{section}\makeatother

\title{\bf QUANTUM MATRIX BALL: THE BERGMAN KERNEL}

\author{\sl D. Shklyarov \and \sl S. Sinel'shchikov \and \sl L. Vaksman}

\date{\tt Institute for Low Temperature Physics \& Engineering \\
National Academy of Sciences of Ukraine}

\newtheorem{theorem}{Theorem}[section]
\newtheorem{lemma}[theorem]{Lemma}
\newtheorem{proposition}[theorem]{Proposition}
\newtheorem{corollary}[theorem]{Corollary}

\begin{document}
\large
\maketitle

\bigskip

\section{Introduction}

 A study of q-analogues for bounded symmetric domains, in particular, for
quantum matrix balls, was started in \cite{SV, SSV}. This work continues
studying these 'balls' and presents the associated Bergman kernels.

 Everywhere below $q \in(0,1)$, $m \le n$, $N=m+n$.

 We assume knowledge of the results of \cite{SSV} and keep the notation of
that work.

\bigskip

\section{\boldmath${\rm Pol}({\rm Mat}_{mn})_{q,y}\simeq{\rm Pol}(X)_{q,x}$}

 The $*$-algebra ${\rm Pol}({\rm Mat}_{mn})_q$, a quantum analogue of
polynomial algebra on the space of matrix, was described in \cite{SSV}. This
algebra was defined in terms of the generators $z_a^\alpha$,
$\alpha=1,\dots,m$; $a=1,\dots,n$, and the commutation relations
\begin{equation}\label{zaa1}z_a^\alpha z_b^\beta-qz_b^\beta
z_a^\alpha=0,\qquad a=b \quad \& \quad \alpha<\beta,\qquad{\rm or}\qquad a<b
\quad \& \quad \alpha=\beta,
\end{equation}
\begin{equation}\label{zaa2}z_a^\alpha z_b^\beta-z_b^\beta
z_a^\alpha=0,\qquad \alpha<\beta \quad \& \quad a>b,
\end{equation}
\begin{equation}\label{zaa3}z_a^\alpha z_b^\beta-z_b^\beta
z_a^\alpha-(q-q^{-1})z_a^\beta z_b^\alpha=0,\qquad \alpha<\beta \quad \&
\quad a<b,
\end{equation}
\begin{equation}(z_b^\beta)^*z_a^\alpha=q^2 \cdot \sum_{a',b'=1}^n
\sum_{\alpha',\beta'=1}^mR_{ba}^{b'a'}R_{\beta \alpha}^{\beta' \alpha'}\cdot
z_{a'}^{\alpha'}(z_{b'}^{\beta'})^*+(1-q^2)\delta_{ab}\delta^{\alpha \beta},
\end{equation}
with $\delta_{ab}$, $\delta^{\alpha \beta}$ being the Kronecker symbols, and
$$R_{ij}^{kl}=\left \{\begin{array}{ccl}q^{-1} &,& i \ne j \quad \& \quad
i=k \quad \& \quad j=l \\ 1 &,& i=j=k=l \\ -(q^{-2}-1) &,& i=j \quad \&
\quad k=l \quad \& \quad l>j \\ 0 &,& {\rm otherwise}\end{array}\right.$$

 The subalgebra generated by $(z_a^\alpha)^*$, $\alpha=1,\dots,m$;
$a=1,\dots,n$, is denoted by ${\mathbb C}[\overline{\rm Mat}_{mn}]_q$, and
the subalgebra generated by $z_a^\alpha$, $\alpha=1,\dots,m$; $a=1,\dots,n$,
is denoted by ${\mathbb C}[{\rm Mat}_{mn}]_q$. Obviously, ${\rm Pol}({\rm
Mat}_{mn})_q={\mathbb C}[{\rm Mat}_{mn}]_q{\mathbb C}[\overline{\rm
Mat}_{mn}]_q$.

 Our work \cite{SSV} provides an embedding of the $*$-algebra ${\rm
Pol}({\rm Mat}_{mn})_q$ into a $*$-algebra of functions on a quantum
principal homogeneous space. Remind the definition of that $*$-algebra.

 Consider the well known algebra ${\mathbb C}[SL_N]_q$ of regular functions
on the quantum group $SL_N$. Its generators are $t_{ij}$, $i,j=1,\dots,N$,
and the complete list of relations includes the relations similar to
(\ref{zaa1}) -- (\ref{zaa3}) and the equality $\det_qT=1$. (Here $\det_qT$
is a q-determinant of the matrix $T=(t_{ij})_{i,j=1,\dots,N}$:
\begin{equation}\label{qdet}\det \nolimits_qT \stackrel{\rm def}{=}\sum_{s
\in S_N}(-q)^{l(s)}t_{1s(1)}\cdot t_{2s(2)}\cdots t_{Ns(N)},
\end{equation}
with $l(s)={\rm card}\{(i,j)|\;i<j \quad \& \quad s(i)>s(j)\}$.)

 Note that q-minors of $T$ are defined similarly to (\ref{qdet}):
$$t_{IJ}^{\wedge k}\stackrel{\rm def}{=}\sum_{s \in
S_k}(-q)^{l(s)}t_{i_1j_{s(1)}}\cdot t_{i_2j_{s(2)}}\cdots t_{i_kj_{s(k)}},$$
with $I=\{(i_1,i_2,\dots,i_k)|\;1 \le i_1<i_2<\dots<i_k \le N \}$,
$J=\{(j_1,j_2,\dots,j_k)|\;1 \le j_1<j_2<\dots<j_k \le N \}$.

 In \cite{SSV} the two $*$-algebras ${\rm Pol}(\widetilde{X})_q=({\mathbb
C}[SL_N]_q,*)$ and ${\mathbb C}[SU_N]_q=({\mathbb C}[SL_N]_q,\star)$ have
been considered, with the involutions $*$ and $\star$ being given by
\begin{equation}\label{invs}{t_{ij}^\star=
(-q)^{j-i}t_{\{1,\dots,\widehat{i},\dots,N \}\{1,\dots,\widehat{j},\dots,N
\}}\atop t_{ij}^*={\rm sign}((i-m+{1 \over 2})(n-j+{1 \over
2}))t_{ij}^\star.}
\end{equation}
They are called the algebra of polynomial functions on a quantum principal
homogeneous space of $SU_{nm}$ and the algebra of regular functions on the
quantum group $SU_N$, respectively. (The latter $*$-algebra is well known
(see \cite{CP})).

 We follow \cite{SSV} in introducing the notation $t=t_{\{1,2,\dots,m
\}\{n+1,n+2,\dots,N \}}$, $x=tt^*$.

 The following lemma is deducible from (\ref{invs}) and a general formula of
Ya. Soibelman \cite{S}, \cite[p. 432]{CP} for the involution $\star$.

\medskip

\begin{lemma} Let ${\rm card}(J)=m$, $J^c=\{1,2,\dots,N \}\setminus J$,
$l(J,J^c)={\rm card}\{(j',j'')\in J \times J^c|\;j'>j'' \}$. Then
\begin{equation}\label{min*}\left(t_{\{1,2,\dots,m \}J}^{\wedge
m}\right)^*=(-1)^{{\rm card}(\{1,2,\dots,n \}\cap
J)}(-q)^{l(J,J^c)}t_{\{m+1,m+2,\dots,N \}J^c}^{\wedge n}.
\end{equation}
\end{lemma}

\medskip

\begin{corollary}\label{tt*}$tt^*=t^*t$.\end{corollary}

\medskip

 Consider a localization ${\rm Pol}(\widetilde{X})_{q,x}$ of the $*$-algebra
${\rm Pol}(\widetilde{X})_q$ with respect to the multiplicative system
$x,x^2,x^3,\dots$ . (The $*$-algebra ${\rm Pol}(\widetilde{X})_{q,x}$ has
no zero divisors. It is derivable from ${\rm Pol}(\widetilde{X})_q$ via
adding a selfadjoint element $x^{-1}$: $xx^{-1}=x^{-1}x=1$,
$(x^{-1})^*=x^{-1}$.)

 A localization ${\mathbb C}[SL_N]_{q,t}$ of the algebra ${\mathbb
C}[SL_N]_q$ with respect to the multiplicative system $t,t^2,\dots$ is
defined in a similar way. Evidently, ${\mathbb
C}[SL_N]_{q,t}\hookrightarrow{\rm Pol}(\widetilde{X})_{q,x}$ since
$t^{-1}=t^*x^{-1}=x^{-1}t^*$ by corollary \ref{tt*}.

 The embedding of $*$-algebras ${\EuScript I}:{\rm Pol}({\rm Mat}_{mn})_q
\hookrightarrow{\rm Pol}(\widetilde{X})_{q,x}$ mentioned above is given by
${\EuScript I}:z_a^\alpha \mapsto t^{-1}t_{\{1,2,\dots,m \}J_{a
\alpha}}^{\wedge m}$, with $J_{a \alpha}=\{n+1,n+2,\dots,N \}\setminus
\{N+1-\alpha \}\cup \{a \}$.

 A crucial point in what follows will be such element $y \in{\rm Pol}({\rm
Mat}_{mn})_q$ that ${\EuScript I}y=x^{-1}$. Our immediate intention is to
construct this $y$.

 We need a notation for q-minors of the matrix
$(z_a^\alpha)_{\alpha=1,\dots,m;\;a=1,\dots,n}$. Suppose $1 \le
\alpha_1<\alpha_2<\dots<\alpha_k \le m$, $1 \le a_1<a_2<\dots<a_k \le n$.
Set up
$$z_{\quad \{a_1,a_2,\dots,a_k \}}^{\wedge k
\{\alpha_1,\alpha_2,\dots,\alpha_k \}}\stackrel{\rm def}{=}\sum_{s \in
S_k}(-q)^{l(s)}z_{a_1}^{\alpha_{s(1)}}z_{a_2}^{\alpha_{s(2)}}\cdots
z_{a_k}^{\alpha_{s(k)}}.$$

 In the proof of the next statement we use the fact that ${\mathbb
C}[SL_N]_q$ and ${\mathbb C}[SL_N]_{q,t}$ are $U_q \mathfrak{sl}_N^{\rm
op}\otimes U_q \mathfrak{sl}_N$-module algebras, ${\mathbb C}[{\rm
Mat}_{mn}]_q$ is a $U_q \mathfrak{sl}_N$-module algebra, and the embedding
${\EuScript I}$ is a morphism of $U_q \mathfrak{sl}_N$-modules. When
considering $U_q \mathfrak{sl}_N$- and $U_q \mathfrak{sl}_N^{\rm
op}$-modules, we use, together with the elements $K_i$, $i=1,\dots,N-1$, the
operators $H_i$, $i=1,\dots,N-1$, introduced in \cite{SSV}. (The
relationship between those is $K_i=q^{H_i}$, $i=1,\dots,N-1$.)

\medskip

\begin{lemma} Let $k \in{\mathbb N}$. There exists $c(q,k)$ such that for
all $1 \le \alpha_1<\alpha_2<\dots<\alpha_k \le m$, $1 \le
a_1<a_2<\dots<a_k \le n$, in the algebra ${\mathbb C}[SL_N]_{q,t}$, one has
\begin{equation}\label{Iz}{\EuScript I}z_{\quad \{a_1,a_2,\dots,a_k
\}}^{\wedge k \{m+1-\alpha_k,m+1-\alpha_{k-1},\dots,m+1-\alpha_1
\}}=c(q,k)t^{-1}t_{\{1,2,\dots,m \}J}^{\wedge m},
\end{equation}
with $J=\{n+1,n+2,\dots,N \}\setminus
\{n+\alpha_1,n+\alpha_2,\dots,n+\alpha_k \}\cup \{a_1,a_2,\dots,a_k \}$.
\end{lemma}

\smallskip

 {\bf Proof.} Consider the linear span of all $z_{\quad \{a_1,a_2,\dots,a_k
\}}^{\wedge k \{m+1-\alpha_k,m+1-\alpha_{k-1},\dots,m+1-\alpha_1 \}}$ and
the linear span of all $t^{-1}t_{\{1,2,\dots,m \}J}^{\wedge m}$, with $k$
being fixed. They are both free modules over the subalgebra $U_q
\mathfrak{sl}_n \otimes U_q \mathfrak{sl}_m \subset U_q \mathfrak{sl}_N$.
Since ${\EuScript I}$ is a morphism of $ U_q \mathfrak{sl}_N$-modules, it
suffices to prove (\ref{Iz}) in the special case as follows:
$${\EuScript I}z_{\quad \{n-k+1,n-k+2,\dots,n \}}^{\wedge k
\{m-k+1,m-k+2,\dots,m \}}=c(q,k)t^{-1}t_{\{1,2,\dots,m \}J_k}^{\wedge m},$$
with $J_k=\{n-k+1,n-k+2,\dots,n,n+k+1,n+k+2,\dots,N \}$.

 Let ${\mathbb F}\subset{\mathbb C}[SL_N]_q$ be a subalgebra generated by
$\{t_{ij}|\;1 \le i \le m,\;1 \le j \le N \}$ and ${\mathbb F}_t
\subset{\mathbb C}[SL_N]_{q,t}$ a subalgebra generated by the same elements
as above and $t^{-1}$.

 It is easy to prove that $f={\EuScript I}z_{\quad \{n-k+1,n-k+2,\dots,n
\}}^{\wedge k \{m-k+1,m-k+2,\dots,m \}}$ belongs to ${\mathbb F}_t$ and is
a solution of the following system of homogeneous linear equations:
\begin{equation}\label{linv}(E_i \otimes 1)f=(F_i \otimes 1)f=(H_i \otimes
1)f=0,\qquad i=1,2,\dots,m-1,
\end{equation}
$$(H_m \otimes 1)f=0,$$
$$(1 \otimes F_j)f=0,\qquad j \ne n,$$
$$(1 \otimes H_j)f=\left \{\begin{array}{ccl}-f &,& j=n \pm k \\ 0 &,& j
\notin \{n,n-k,n+k \}\end{array}\right..$$
Since $t^{-1}t_{\{1,2,\dots,m \}J_k}^{\wedge m}$ satisfies all the above
equations, it suffices to prove that the space of solutions $f \in{\mathbb
F}_t$ of this system is one-dimensional. Arguing just as in the proof of
\cite[lemma 8.2]{SSV}, we obtain the following results. The subalgebra of
solutions $f \in{\mathbb F}_t$ of (\ref{linv}) is generated by $t^{-1}$,
$t_{\{1,2,\dots,m \}I}^{\wedge m}$, ${\rm card}\,I=m$. The subalgebra
${\mathbb F}_{\rm inv}$ of solutions of (\ref{linv}) together with $(H_m
\otimes 1)f=0$, is generated by ratios of quantum minors
$t^{-1}t_{\{1,2,\dots,m \}I}^{\wedge m}$, ${\rm card}\,I=m$, and the
subalgebra ${\mathbb F}_{\rm prim}=\{f \in{\mathbb F}_{\rm inv}|\;(1 \otimes
F_j)f=0,\:j \ne m \}$ is generated by ratios of quantum minors
$t^{-1}t_{\{1,2,\dots,m \}J_l}^{\wedge m}$,
$J_l=\{n-l+1,n-l+2,\dots,n,n+l+1,n+l+2,\dots,N \}$, $l \le m$.

 What remains is to consider the linear span ${\mathfrak h}'$ of $1 \otimes
H_j$, $j \ne n$, and to elaborate the linear independence in $({\mathfrak
h}')^*$ of the weights of generators in ${\mathbb F}_{\rm prim}$:
$$(1 \otimes H_j)\left(t^{-1}t_{\{1,2,\dots,m \}J_l}^{\wedge
m}\right)=\left \{\begin{array}{ccl}-t^{-1}t_{\{1,2,\dots,m \}J_l}^{\wedge
m} &,& j=n \pm l \\ 0 &,& j \notin \{n,n-l,n+l \}\end{array}\right..\eqno
\blacksquare$$

\medskip

 The action of the subalgebra $U_q \mathfrak{sl}_N \hookrightarrow U_q
\mathfrak{sl}_N^{\rm op}\otimes U_q \mathfrak{sl}_N$ is to be referred to in
the sequel more extensively than that of the subalgebra $U_q
\mathfrak{sl}_N^{\rm op} \hookrightarrow U_q \mathfrak{sl}_N^{\rm op}\otimes
U_q \mathfrak{sl}_N$. Just as in \cite{SSV}, we write $\xi f$ instead of $(1
\otimes \xi)f$ in all the cases where this could not lead to a confusion.

 Apply the relations $F_nz_a^\alpha=q^{1/2}\delta_{an}\delta_{\alpha m}$,
$H_nz_b^\beta=0$ for $b \ne n$ and $\beta \ne m$ (see \cite{SSV}) to get
$$F_nz_{\quad \{n-k+1,\dots,n \}}^{\wedge k \{m-k+1,\dots,m
\}}=q^{1/2}z_{\quad \{n-k+1,\dots,n-1 \}}^{\wedge k-1 \{m-k+1,\dots,m-1
\}}.$$
It follows from $\Delta(F_n)=F_n \otimes K_n^{-1}+1 \otimes F_n$,
$F_n(t^{-1})=0$ that
$$F_n \left(t^{-1}t_{\{1,2,\dots,m \}\{n-k+1,\dots,n,n+k+1,\dots,N
\}}^{\wedge m}\right)=t^{-1}F_nt_{\{1,2,\dots,m
\}\{n-k+1,\dots,n,n+k+1,\dots,N \}}^{\wedge m}=$$
$$=q^{1/2}t^{-1}t_{\{1,2,\dots,m \}\{n-k+1,\dots,n-1,n+1,n+k+1,\dots,N
\}}^{\wedge m}.$$
Hence, $c(q,k)=c(q,k-1)=\dots=c(q,1)=1$, and thus we have proved

\medskip

\begin{proposition}\label{Izm}Let $1 \le \alpha_1<\alpha_2<\dots<\alpha_k
\le m$, $1 \le a_1<a_2<\dots<a_k \le n$, $J=\{n+1,n+2,\dots,N \}\setminus
\{n+\alpha_1,n+\alpha_2,\dots,n+\alpha_k \}\cup \{a_1,a_2,\dots,a_k \}$.
Then
$${\EuScript I}z_{\quad \{a_1,a_2,\dots,a_k \}}^{\wedge k
\{m+1-\alpha_k,m+1-\alpha_{k-1},\dots,m+1-\alpha_1
\}}=t^{-1}t_{\{1,2,\dots,m \}J}^{\wedge m}$$
\end{proposition}

\medskip

 Turn back to a construction of such $y \in{\rm Pol}({\rm Mat}_{mn})_q$ that
${\EuScript I}y=x^{-1}$. It follows from (\ref{qdet}), (\ref{min*}) that
\begin{equation}\label{tmtm*}\sum_{J \subset \{1,\dots,N \}\atop{\rm
card}(J)=m}(-1)^{{\rm card}(\{1,2,\dots,n \}\cap J)}t_{\{1,2,\dots,m
\}J}^{\wedge m}\left(t_{\{1,2,\dots,m \}J}^{\wedge m}\right)^*=1.
\end{equation}

 The following is due to proposition \ref{Izm}, (\ref{tmtm*}), and the
injectivity of ${\EuScript I}$.

\medskip

\begin{theorem}\label{yexp} There exists a unique element $y \in{\rm
Pol}({\rm Mat}_{mn})_q$ such that ${\EuScript I}y=x^{-1}$. It is given
explicitly by
$$y=1+\sum_{k=1}^m(-1)^k \sum_{\{J'|\;{\rm card}(J')=k \}}\sum_{\{J''|\;{\rm
card}(J'')=k \}}z_{\quad J''}^{\wedge kJ'}\left(z_{\quad J''}^{\wedge
kJ'}\right)^*,$$
with $J'\subset \{1,2,\dots,m \}$, $J''\subset \{1,2,\dots,n \}$.
\end{theorem}

\medskip

{\sc Example 2.6.} In the case of quantum ball in ${\mathbb C}^n$ $(m=1)$,
one has $z_a=z_a^1$, $y=1-\sum \limits_{a=1}^nz_az_a^*$.

\medskip \stepcounter{theorem}

\begin{corollary}\label{zy} For all $\alpha=1,\dots,m$, $a=1,\dots,n$, one
has $$z_a^\alpha y=q^{-2}yz_a^\alpha,\qquad
(z_a^\alpha)^*y=q^2y(z_a^\alpha)^*.$$
\end{corollary}

\medskip

 Note that $f \in{\rm Pol}({\rm Mat}_{mn})_q$, $f \ne 0$ imply $yf \ne 0$,
$fy \ne 0$ because of the injectivity of the map ${\EuScript I}:{\rm
Pol}({\rm Mat}_{mn})_q \to{\rm Pol}(\widetilde{X})_{q,x}$. Hence ${\rm
Pol}({\rm Mat}_{mn})_q \hookrightarrow{\rm Pol}({\rm Mat}_{mn})_{q,y}$.

 The work \cite{SSV} introduces a subalgebra ${\rm Pol}(X)_q$ of all $U_q
\mathfrak{s}(\mathfrak{gl}_m \times \mathfrak{gl}_n)^{\rm op}$-invariants of
${\rm Pol}(\widetilde{X})_q$. Let ${\rm Pol}(X)_{q,x}$ and ${\rm Pol}({\rm
Mat}_{mn})_{q,y}$ be respectively localizations of the $*$-algebras ${\rm
Pol}(X)_q$ and ${\rm Pol}({\rm Mat}_{mn})_q$ with respect to the
multiplicative systems $x^{\mathbb N}$, $y^{\mathbb N}$. Now theorem
\ref{yexp} and \cite[proposition 8.1]{SSV} provide a 'canonical'
isomorphism

\medskip

\begin{proposition} ${\rm Pol}({\rm Mat}_{mn})_{q,y}
{\to_{\!\!\!\!\!\!\!\!\!_{\textstyle \sim \atop \scriptstyle{\,{\EuScript
I}\,}}}}{\rm Pol}(X)_{q,x}$.
\end{proposition}

\medskip

 (The injectivity of ${\EuScript I}$ is evident. In fact, it follows from
${\EuScript I}\left(\sum \limits_{k=0}^My^{-k}f_k \right)=0$,
$f_0,f_1,\dots,f_M \in{\rm Pol}({\rm Mat}_{mn})_q$, that ${\EuScript
I}\left(\sum \limits_{k=0}^My^{M-k}f_k \right)=0$, and so $\sum
\limits_{k=0}^My^{M-k}f_k=0$, $\sum \limits_{k=0}^My^{-k}f_k=0$.)

\medskip

 {\sc Remark 2.9.} In some contexts the generators $z_a^\alpha$ become
inconvenient and are to be replaced by $z_{\alpha
a}=(-q)^{\alpha-1}z_a^{m+1-\alpha}$, $\alpha=1,\dots,m$, $a=1,\dots,n$.
(This passage could treated as drawing down the 'Greek' index via the
'tensor' $\varepsilon=(-q)^{\alpha-1}\delta_{\alpha+\beta,m+1}$.)

\medskip \stepcounter{theorem}

\begin{proposition} Consider the matrix $Z=(z_{\alpha
a})_{\alpha=1,\dots,m,a=1,\dots,n}$. In the matrix algebra with entries from
${\mathbb C}[SL_N]_{q,t}$ one has
\begin{equation}\label{IZ}{\EuScript I}(Z)=T_{12}^{-1}T_{11},\end{equation}
with ${\EuScript I}(Z)=({\EuScript I}(z_{\alpha a}))$, $T_{11}=(t_{\alpha
a})$, $T_{12}=(t_{\alpha n+\beta})$, $\alpha,\beta=1,\dots,m,a=1,2,\dots,n$.
\end{proposition}

\smallskip

 {\bf Proof.} Let $T=(t_{ij})_{i,j=1,\dots,m}$ and $\det_q'T=\sum_{s \in
S_m}(-q)^{-l(s)}t_{s(m)m}t_{s(m-1)m-1}\cdots t_{s(1)1}$. We are about to
prove that in ${\mathbb C}[{\rm Mat}_{mn}]_q$ one has
\begin{equation}\label{det'}\det \nolimits_q'T=\det
\nolimits_qT.\end{equation}
The relation $\det_q'T={\tt const}(q)\det_qT$ follows from the $U_q
\mathfrak{sl}_m$-invariance of $\det_q'T$, $\det_qT$, if one also takes
into account that the spaces of homogeneous degree $m$ invariants are
one-dimensional. To compute the constant ${\tt const}(q)$, it suffices to
pass to the quotient algebra with respect to the bilateral ideal generated
by $t_{ij}$, $i \ne j$.

 Now (\ref{IZ}) is derivable from (\ref{det'}) and the explicit form of
$T_{12}^{-1}$:
$$(T_{12}^{-1})_{\alpha \beta}=(\det
\nolimits_qT_{12})^{-1}\cdot(-q)^{\alpha-\beta}\det
\nolimits_q((T_{12})_{\beta \alpha}),\qquad \alpha,\beta=1,\dots,m$$
(Here, just as in the classical case $q=1$, $(T_{12})_{\beta \alpha}$ is a
matrix derived from $T_{12}$ by discarding the line $\beta$ and column
$\alpha$.) \hfill $\blacksquare$

\medskip

 {\sc Remark 2.11.} In the classical limit ($q=1$) one has
$$y=1+\sum_{k=1}^m(-1)^k{\rm tr}(Z^{\wedge k}(Z^*)^{\wedge
k})=1+\sum_{k=1}^m(-1)^k{\rm tr}((ZZ^*)^{\wedge k})=\det(1-ZZ^*).$$
(These relations are evident since their proof reduces to considering the
special case $z_{\alpha a}=\lambda_ \alpha \delta_{\alpha a}$,
$\lambda_ \alpha \in{\mathbb C}$.)

\bigskip

\section{Hardy-Bergman spaces}

 Remind some results of \cite{SSV}. An extension ${\rm Fun}({\mathbb U})_q$
of the covariant $*$-algebra ${\rm Pol}({\rm Mat}_{mn})_q$ was produced
there via adding to the list of its generators $\{z_a^\alpha \}$ of such
element $f_0$ that $f_0=f_0^2=f_0^*$, and
$(z_a^\alpha)^*f_0=f_0z_a^\alpha=0$ for all $\alpha=1,\dots,m$,
$a=1,\dots,n$.

 Let $D({\mathbb U})_q={\rm Fun}({\mathbb U})_qf_0{\rm Fun}({\mathbb U})_q$
be the $*$-algebra of finite functions in the quantum ball, and ${\cal
H}={\rm Fun}({\mathbb U})_qf_0={\mathbb C}[{\rm Mat}_{mn}]_qf_0$. It follows
from the explicit formulae for the action of $U_q \mathfrak{sl}_N$ in ${\rm
Fun}({\mathbb U})_q$ (see \cite{SSV}) that $D({\mathbb U})_q$ is a $U_q
\mathfrak{su}_{nm}$-module algebra, and ${\cal H}$ is a $U_q
\mathfrak{b}_+$-module algebra. (Here $U_q \mathfrak{b}_+\subset U_q
\mathfrak{sl}_N$ is a subalgebra generated by $K_j^{\pm 1}$, $E_j$,
$j=1,\dots,N-1$). We follows \cite{SSV} in denoting the 'natural'
representations on ${\rm Fun}({\mathbb U})_q$ and $U_q \mathfrak{b}_+$ in
${\cal H}$ by $\Theta$ and $\Gamma$ respectively. In section 9 of that work
an explicit formula for the invariant integral was obtained:
\begin{equation}\label{ii}\int \limits_{{\mathbb U}_q}f d \nu={\rm
tr}(\Theta(f)\Gamma(e^{h \check{\rho}})),\qquad f \in D({\mathbb U})_q,
\end{equation}
with $\check{\rho}=\dfrac{1}{2}\sum \limits_{j=1}^{N-1}j(N-j)H_j$,
$h>0$, and the operators $\Gamma(H_j)$ in ${\cal H}$ are determined by
$q=e^{-h/2}$, $\Gamma(K_j^{\pm 1})=q^{\pm \Gamma(H_j)}$, $j,\dots,N-1$. (We
follow an agreement of \cite{SSV} in restricting to considering only those
$U_q \mathfrak{sl}_N$-modules which admit well defined actions of the
'elements' $H_j$, $X_j^+=E_jq^{-{1 \over 2}H_j}$, $X_j^-=q^{{1 \over
2}H_j}F_j$, $j=1,\dots,N-1$.)

 Among the main results of \cite{SSV} one should mention, in particular, the
positivity of the scalar product in ${\cal H}$ defined by
\begin{equation}\label{sp}(\psi_1,\psi_2)f_0=\psi_2^*\psi_1,\qquad
\psi_1,\psi_2 \in{\cal H}.
 \end{equation}
The positivity of the scalar product (\ref{sp}) implies the positivity  of
the invariant integral (\ref{ii}) (see \cite{SSV}).

 Remind (see \cite{SSV}) that ${\mathbb C}[{\rm Mat}_{mn}]_q=\bigoplus
\limits_{k=0}^\infty{\mathbb C}[{\rm Mat}_{mn}]_{q,k}$, ${\mathbb C}[{\rm
Mat}_{mn}]_{q,k}=\{f|\;\deg \,f=k \}$, ${\cal H}=\bigoplus
\limits_{k=0}^\infty{\cal H}_k$, ${\cal H}_k={\mathbb C}[{\rm
Mat}_{mn}]_{q,k}f_0$. The equality $yf_0=f_0$ and the commutation relations
$yz_a^\alpha=q^2z_a^\alpha y$, $\alpha=1,\dots,m$, $a=1,\dots,n$ imply

\medskip

\begin{lemma} For all $k,\lambda \in{\mathbb Z}_+$ one has
\begin{equation}\label{Thy}\left.\Theta(y)^\lambda \right|_{{\cal
H}_k}=q^{2k \lambda}I.
\end{equation}
\end{lemma}

\medskip

 The non-integral powers $\Theta(y)^\lambda \in{\rm End}({\cal H})$ of the
operator $\Theta(y)$ are defined by (\ref{Thy}).

\begin{lemma}\label{cl} For $\lambda>N-1$ one has
$${\rm tr}\left(\Theta(y)^\lambda \Gamma(e^{h
\check{\rho}})\right)=\prod_{j=0}^{n-1}\prod_{k=0}^{m-1}
\left(1-q^{2(\lambda+1-N)}q^{2(j+k)}\right)^{-1}$$
\end{lemma}

\smallskip

 {\bf Proof.} Equip the linear span $\mathfrak h$ of the 'elements' $H_j$,
$j=1,2,\dots,N-1$, with a scalar product
$$(H_i,H_j)=\left \{\begin{array}{ccl}2 &,& i-j=0 \\ -1 &,& |i-j|=1 \\ 0
&,& |i-j|>1 \end{array}\right..$$

 Let $\{\alpha_j \}$ be the standard basis of simple roots in ${\mathfrak
h}^*$. Then $\alpha_j(\check{\rho})=(H_j,\check{\rho})=1$ for all
$j=1,\dots,N-1$, and the relation
$$\check{\rho}z_n^m={1 \over
2}\left((m+1)(n-1)H_{n-1}+(m-1)(n+1)H_{n+1}+mnH_n \right)z_n^m=$$
$$={1 \over 2}\left(-(m+1)(n-1)-(m-1)(n+1)+2mn \right)z_n^m=z_n^m$$
implies
\begin{equation}\label{roz}\check{\rho}z_a^\alpha=
(N+1-a-\alpha)z_a^\alpha,\qquad \alpha=1,\dots,m,\;a=1,\dots,n.
\end{equation}

 The vectors $\{(z_1^1)^{k_{11}}\cdots(z_n^m)^{k_{mn}}f_0 \}$, $k_{a
\alpha}\in{\mathbb Z}_+$, $\alpha=1,\dots,m$, $a=1,\dots,n$, form a basis in
${\cal H}$, and, by a virtue of (\ref{roz}), one has
\begin{equation}\label{Gam}\Gamma(e^{h
\check{\rho}})((z_1^1)^{k_{11}}\cdots(z_n^m)^{k_{mn}}f_0)=q^{-2 \sum k_{a
\alpha}(N+1-a-\alpha)}((z_1^1)^{k_{11}}\cdots(z_n^m)^{k_{mn}}f_0).
\end{equation}

 What remains is to apply the definition of the operator
$\Theta(y)^\lambda$:
\begin{equation}\label{fe}{\rm tr}\left(\Theta(y)^\lambda \Gamma(e^{h
\check{\rho}})\right)=\sum_{k_{11}=0}^\infty \cdots \sum_{k_{mn}=0}^\infty
q^{2 \sum \limits_{a=1}^n \sum
\limits_{\alpha=1}^m(\lambda-(N+1-a-\alpha))k_{a \alpha}}
\end{equation}
and to sum the geometric progressions in (\ref{fe}). \hfill $\blacksquare$

\medskip

\begin{proposition}\label{btf} For any $f \in{\rm Fun}({\mathbb U})_q$,
$\Theta(f)$ is a bounded operator in the pre-Hilbert space ${\cal H}$.
\end{proposition}

\smallskip

 {\bf Proof.} By a virtue of \cite[remark 8.6]{SSV}, it suffices to consider
the case $f \in{\rm Pol}({\rm Mat}_{mn})_q$. In the work alluded above, a
$*$-representation $\Pi$ was constructed; it was also shown to be unitary
equivalent to the $*$-representation $\Theta$. Thus, the inequality
$\|\Theta(f)\|<\infty$ follows from $\|\Pi(f)\|<\infty$ (see the appendix).
\hfill $\blacksquare$

\medskip

 By proposition \ref{btf}, for $\lambda>N-1$ one has a well defined linear
functional
$$\int \limits_{{\mathbb U}_q}fd \nu_ \lambda \stackrel{\rm
def}{=}C(\lambda){\rm tr}\left(\Theta(f)\Theta(y)^\lambda \Gamma(e^{h
\check{\rho}})\right),\qquad f \in{\rm Fun}({\mathbb U})_q,$$
with
\begin{equation}\label{Cl}C(\lambda)=\prod_{j=0}^{n-1}\prod_{k=0}^{m-1}
\left(1-q^{2(\lambda+1-N)}q^{2(j+k)}\right).
\end{equation}

\medskip

\begin{proposition}\label{pnl} For all $\lambda>N-1$, the linear functional
$\int \limits_{{\mathbb U}_q}fd \nu_ \lambda$ on the $*$-algebra ${\rm
Fun}({\mathbb U})_q$ is positive, and $\int \limits_{{\mathbb U}_q}1d \nu_
\lambda=1$.
\end{proposition}

\smallskip

 {\bf Proof.} The latter relation follows from lemma \ref{cl}. To verify the
positivity of $\int \limits_{{\mathbb U}_q}fd \nu_ \lambda$, it suffices to
observe that the $*$-representation $\Theta$ of ${\rm Fun}({\mathbb U})_q$
is faithful, and the bounded operator $C(\lambda)\Theta(y)^\lambda \Gamma
\left(e^{h \check{\rho}}\right)=C(\lambda)\Gamma \left(e^{h
\check{\rho}}\right)\Theta(y)^\lambda$ is positive. \hfill $\blacksquare$

\medskip

 Consider a completion $L^2(d \nu_ \lambda)_q$ of the vector space ${\rm
Fun}({\mathbb U})_q$ with respect to the norm $\|f \|_ \lambda=\left(\int
\limits_{{\mathbb U}_q}f^*fd \nu_ \lambda \right)^{1/2}$. The closure
$L_a^2(d \nu_ \lambda)_q$ of the linear subvariety ${\mathbb C}[{\rm
Mat}_{mn}]_q$ in the Hilbert space $L^2(d \nu_ \lambda)_q$ will be called
the Hardy-Bergman space.

 There exists a very useful approach in which the $*$-algebra ${\rm
Pol}({\rm Mat}_{mn})_q$ is treated as a q-analogue of the Weil algebra, and
the operators $(1-q^2)^{-1/2}T(z_a^\alpha)$, $(1-q^2)^{-1/2}T(z_a^\alpha)^*$
as q-analogues of creation and annihilation operators respectively. In this
context, the following result becomes a q-analogue of the Stone-von-Neumann
theorem.

\medskip

\begin{theorem} There exists a faithful irreducible $*$-representation of
${\rm Pol}({\rm Mat}_{mn})_q$ by bounded operators in a Hilbert space. This
representation is unique up to unitary equivalence.
\end{theorem}

\smallskip

 {\bf Proof.} By proposition \ref{btf}, there exists a well defined
$*$-representation $\overline{\Theta}$ of ${\rm Pol}({\rm Mat}_{mn})_q$ in
a completion $\overline{\cal H}$ of the pre-Hilbert space ${\cal H}$. This
$*$-representation is faithful by \cite[proposition 8.8]{SSV}. Furthermore,
if a bounded linear operator $A$ commutes with all the operators
$\overline{\Theta}(f)$, $f \in{\rm Pol}({\rm Mat}_{mn})_q$, then, in
particular $\overline{\Theta}(y)A=A \overline{\Theta}(y)$. Hence $Af_0=af_0$
for some $a \in{\mathbb C}$ since ${\mathbb C}f_0$ is an eigenspace of
$\overline{\Theta}(y)$. It follows that $A=aI$. That is, $\overline{\Theta}$
is irreducible. What remains is to demonstrate a uniqueness of the faithful
irreducible $*$-representation.

 Let $T$ be a faithful irreducible $*$-representation of ${\rm Pol}({\rm
Mat}_{mn})_q$ by bounded linear operators in a Hilbert space. The same idea
as in \cite{VS} can be used to prove that the non-zero spectrum of the
selfadjoint operator $T(y)$ is discrete. Consider some eigenvector $v$ of
$T(y)$ associated to a largest modulus eigenvalue. By a virtue of corollary
\ref{zy}, $T((z_a^\alpha)^*)v=0$, $\alpha=1,2,\dots,m$, $a=1,2,\dots,n$. It
is easy to show that the kernels of the linear functionals $(T(f)v,v)$,
$(\Theta(f)f_0,f_0)$ on ${\rm Pol}({\rm Mat}_{mn})_q$ are just the same
subspace $\bigoplus \limits_{(j,k)\ne(0,0)}{\mathbb C}[{\rm
Mat}_{mn}]_{q,j}{\mathbb C}[\overline{\rm Mat}_{mn}]_{q,-k}$. Thus
$(T(f)v,v)={\tt const}(\Theta(f)f_0,f_0)$, ${\tt const}>0$, and hence the
map $f_0 \mapsto({\tt const})^{-1/2}v$ admits an extension up to a unitary
map which intertwines the representations $\overline{\Theta}$ and $T$.
\hfill $\blacksquare$

\bigskip

\section{Distributions}

 It will be proved below that the orthogonal projection $P_ \lambda$ in the
Hilbert space $L^2(d \nu_ \lambda)_q$ onto the subspace $L_a^2(d \nu_
\lambda)_q \subset L^2(d \nu_ \lambda)_q$ is an integral operator. Our
principal intention is find the kernel of this integral operator.

 We follow \cite{SSV} with beginning the construction of distributions in
the quantum matrix ball via a completion operation.

 Impose three topologies in ${\rm Pol}({\rm Mat}_{mn})_q$ and prove their
equivalence.

 Associate to each element $\psi \in D({\mathbb U})_q$ a linear functional
$l_ \psi:{\rm Pol}({\rm Mat}_{mn})_q \to{\mathbb C}$ given by $l_
\psi(f)=\int \limits_{{\mathbb U}_q}f \psi d \nu$. Let ${\cal T}$ be the
weakest among the topologies in ${\rm Pol}({\rm Mat}_{mn})_q$ in which all
the linear functionals $l_ \psi$, $\psi \in D({\mathbb U})_q$, are
continuous.

 Consider the {\sl finite dimensional} subspaces ${\rm Pol}({\rm
Mat}_{mn})_{q,i,j}={\mathbb C}[{\rm Mat}_{mn}]_{q,i}\cdot{\mathbb
C}[\overline{\rm Mat}_{mn}]_{q,j}$, $i \le 0$, $j \ge 0$, of the vector
space ${\rm Pol}({\rm Mat}_{mn})_q$, together with such linear operators
$p_{ij}:{\rm Pol}({\rm Mat}_{mn})_q \to{\rm Pol}({\rm Mat}_{mn})_{q,i,j}$
that $f=\sum \limits_{i \ge 0,\;j \le 0}p_{ij}(f)$ for all $f \in{\rm
Pol}({\rm Mat}_{mn})_q$. Let ${\cal T}_1$ be the weakest among the topologies
in ${\rm Pol}({\rm Mat}_{mn})_q$ in which all $p_{ij}$, $i \le 0$, $j \ge
0$, are continuous.

 Consider the {\sl finite dimensional} orthogonal projections $P_i$ in
${\cal H}$ onto the homogeneous components ${\cal H}_i$, $i \in{\mathbb
Z}_+$, and the linear operators $\Theta_{ij}:{\rm Pol}({\rm Mat}_{mn})_q
\to{\rm Hom}({\cal H}_j,{\cal H}_i)$, $f \mapsto P_i \Theta(f)|_{{\cal
H}_j}$, $i,j \in{\mathbb Z}_+$.

 Let ${\cal T}_2$ be the weakest among the topologies in which all the
linear operators $\Theta_{ij}$, $i,j \in{\mathbb Z}_+$, are continuous.

 The completion of ${\rm Pol}({\rm Mat}_{mn})_q$ with respect to ${\cal
T}$ will be denoted by $D({\mathbb U})_q'$ and called the space of
distributions in the quantum matrix ball. The pairing ${\rm Pol}({\rm
Mat}_{mn})_q \times D({\mathbb U})_q \to{\mathbb C}$, $f \times \psi \mapsto
\int \limits_{{\mathbb U}_q}f \psi d \nu$, is extendable up to a pairing
$D({\mathbb U})_q'\times D({\mathbb U})_q \to{\mathbb C}$, $f \times \psi
\mapsto \int \limits_{{\mathbb U}_q}f \psi d \nu$; this justifies the use
of the term 'distribution'.

 One can replace the topology ${\cal T}$ in the definition of $D({\mathbb
U})_q'$ with either ${\cal T}_1$ or ${\cal T}_2$, as it follows from

\medskip

\begin{proposition} The topologies ${\cal T}$, ${\cal T}_1$, ${\cal T}_2$
in ${\rm Pol}({\rm Mat}_{mn})_q$ are equivalent.
\end{proposition}

\smallskip

 {\bf Proof.} By \cite[remark 8.6]{SSV}, $D({\mathbb U})_q
\xrightarrow[\Theta]{\textstyle \sim}{\rm End}({\cal H})_f$, with ${\rm
End}({\cal H})_f \simeq{\cal H}\otimes{\cal H}^*$. Thus, the equivalence of
the topologies ${\cal T}$, ${\cal T}_2$ follows from (\ref{ii}), (\ref{fe}).

 What remains is to prove the equivalence of ${\cal T}_1$ and ${\cal T}_2$.
Consider the linear span ${\cal L}_1$ of the images under the embedding
$({\rm Pol}({\rm Mat}_{mn})_q)_{ij})^*\hookrightarrow({\rm Pol}({\rm
Mat}_{mn})_q)^*$ and the linear span ${\cal L}_2$ of the images under the
embedding $({\rm Hom}({\cal H}_j,{\cal H}_i))^*\hookrightarrow({\rm
Pol}({\rm Mat}_{mn})_q)^*$. It suffices to prove that ${\cal L}_1={\cal
L}_2$. The inclusion ${\cal L}_1 \supseteq{\cal L}_2$ follows from ${\rm
Hom}({\cal H}_j,{\cal H}_i)^*\subseteq \bigoplus
\limits_{k=0}^{\min(i,j)}{\rm Pol}({\rm Mat}_{mn})_{q,i-k,-j+k}^*$. The
converse inclusion follows from ${\rm Pol}({\rm
Mat}_{mn})_{q,i-j}^*\subseteq \bigoplus \limits_{k=0}^{\min(i,j)}{\rm
Hom}({\cal H}_j,{\cal H}_i)^*$. The latter inclusion is easily deducible
from the previous one and \cite[lemma 8.7]{SSV}. \hfill $\blacksquare$

\medskip

 The equivalence of ${\cal T}$ and ${\cal T}_1$ allows one to identify the
topological vector space $D({\mathbb U})_q'$ and the space of formal series
$f=\sum \limits_{i \ge 0,\;j \le 0}f_{ij}$, $f_{ij}\in{\rm Pol}({\rm
Mat}_{mn})_{q,i,j}$ equipped with the topology of coefficientwise
convergence. The structure of a covariant ${\rm Pol}({\rm
Mat}_{mn})_q$-bimodule is transferred by a continuity from ${\rm Pol}({\rm
Mat}_{mn})_q$ onto the above space $D({\mathbb U})_q'$ of formal series.

 Let $\overline{\rm End}({\cal H})=\underset{i,j \ge
0}{\times}{\rm Hom}({\cal H}_j,{\cal H}_i)$ -- the direct product in the
category of vector spaces, i.e. the corresponding space of formal series.
Evidently, ${\rm End}({\cal H})_f \hookrightarrow{\rm End}({\cal
H})\hookrightarrow \overline{{\rm End}({\cal H})}$.

 Consider the embedding $i:D({\mathbb U})_q \hookrightarrow D({\mathbb
U})_q'$ determined via the isomorphisms $D({\mathbb U})_q \simeq{\rm
End}({\cal H})_f$, $D({\mathbb U})_q'\simeq \overline{{\rm End}({\cal H})}$.
(The second isomorphism is a consequence of the equivalence of ${\cal T}$
and ${\cal T}_2$.)

\medskip

\begin{proposition} The embedding of vector spaces $i:D({\mathbb U})_q
\hookrightarrow D({\mathbb U})_q'$ is a morphism of {\bf covariant} ${\rm
Pol}({\rm Mat}_{mn})_q$-modules.
\end{proposition}

\smallskip

 {\bf Proof.} By the construction, $i$ is a morphism of ${\rm Pol}({\rm
Mat}_{mn})_q$-bimodules. What remains is to prove that this is a morphism
of $U_q \mathfrak{sl}_N$-modules. Consider the element $f_0 \in D({\mathbb
U})_q'$ which is the image of $f_0 \in D({\mathbb U})_q$ under the embedding
$i$. It suffices to demonstrate that the relations from \cite{SSV} which
define the structure of $U_q \mathfrak{sl}_N$-module in $D({\mathbb U})_q$,
are also valid in $D({\mathbb U})_q'$:
$$H_nf_0=0,\qquad F_nf_0=\frac{q^{1/2}}{q^{-2}-1}f_0(z_n^m)^*,\qquad
E_nf_0=-\frac{q^{1/2}}{1-q^2}z_n^mf_0,$$
$$H_kf_0=F_kf_0=E_kf_0=0 \quad{\rm for}\quad k \ne n.$$

It is easy to see that
$${\mathbb C}f_0=\{f \in D({\mathbb
U})_q'|\;(z_a^\alpha)^*f=fz_a^\alpha=0,\quad \alpha=1,\dots,m;\;a=1,\dots,n
\}.$$
It follows from the covariance of the ${\rm Pol}({\rm Mat}_{mn})_q$-bimodule
$D({\mathbb U})_q'$ that the subspace ${\mathbb C}f_0$ is a $U_q
\mathfrak{s}(\mathfrak{gl}_n \times \mathfrak{gl}_m)$-submodule. Hence,
$H_nf_0=0$, $H_kf_0=F_kf_0=E_kf_0=0$, $k \ne n$. Similarly,
$${\mathbb C}z_n^mf_0=\{f \in D({\mathbb U})_q'|\;H_0f=2f \;\&\;F_jf=0
\quad{\rm for}\quad j \ne n$$
$$\&\;fz_a^\alpha=0 \quad{\rm for}\quad \alpha=1,\dots,m;\;a=1,\dots,n \}.$$
$${\mathbb C}f_0(z_n^m)^*=\{f \in D({\mathbb U})_q'|\;H_0f=-2f \;\&\;E_jf=0
\quad {\rm for}\quad j \ne n$$
$$\&\;(z_a^\alpha)^*f=0 \quad{\rm for}\quad \alpha=1,\dots,m;\;a=1,\dots,n
\}.$$

 Apply the covariance of the ${\rm Pol}({\rm Mat}_{mn})_q$-bimodule
$D({\mathbb U})_q'$ to get
$$F_nf_0={\tt const}_1f_0(z_n^m)^*,\qquad E_nf_0={\tt const}_2z_n^mf_0.$$
What remains is to prove that ${\tt const}_1=-\dfrac{q^{1/2}}{q^{-2}-1}$,
${\tt const}_2=-\dfrac{q^{1/2}}{1-q^2}$.

 The first constant is accessible from the relations
$$F_n(f_0z_n^m)=F_n0=0,\qquad
f_0(1-(z_n^m)^*z_n^m)=q^2f_0(1-z_n^m(z_n^m)^*)=q^2f_0,$$
and the second one follows from
$$E_n((z_n^m)^*f_0)=E_n0=0,\qquad
(1-(z_n^m)^*z_n^m)f_0=q^2f_0(1-z_n^m(z_n^m)^*)f_0=q^2f_0.$$
(A detailed exposition of these calculations can be found in \cite{SSV}).
\hfill $\blacksquare$

\medskip

 We identify in the sequel finite functions $f \in D({\mathbb U})_q$ with
their images $i(f)\in D({\mathbb U})_q'$ under the embedding $i$.

\medskip \stepcounter{theorem}

{\sc Remark 4.3.} By definition, the linear subspace $D({\mathbb U})_q
\xrightarrow[\Theta]{\textstyle \sim}{\rm End}({\cal H})_f$ is dense in the
topological vector space $D({\mathbb
U})_q'\xrightarrow[\Theta]{\textstyle \sim}\overline{\rm End}({\cal H})$.
The structure of the $D({\mathbb U})_q$-bimodule is extendable by a
continuity from this dense linear subspace onto the entire space $D({\mathbb
U})_q'$. $D({\mathbb U})_q$ is a covariant algebra, hence the $D({\mathbb
U})_q$-bimodule we have obtained is also covariant.

\medskip

 To conclude, we prove the following

\medskip

\begin{proposition}\label{fc} A distribution $f \in D({\mathbb U})_q'$ is a
finite function iff $f{\mathbb C}[{\rm Mat}_{mn}]_{q,M}={\mathbb
C}[\overline{\rm Mat}_{mn}]_{q,-M}f=0$ for some $M \in{\mathbb N}$.
\end{proposition}

\smallskip

 {\bf Proof.} Apply the fact that the linear map ${\mathbb C}[{\rm
Mat}_{mn}]_q \to{\cal H}$, $\psi \mapsto \psi f_0$, is one-to-one, and
${\mathbb C}[{\rm Mat}_{mn}]_q$ is a domain \cite{DD}. For all $k
\in{\mathbb Z}_+$, the linear span of $\Theta(z_a^\alpha){\cal H}_k$,
$\alpha=1,\dots,m$, $a=1,\dots,n$, coincides with ${\cal H}_{k+1}$, and the
linear span of $\Theta((z_a^\alpha)^*){\cal H}_k$, $\alpha=1,\dots,m$,
$a=1,\dots,n$, coincides with ${\cal H}_{k-1}$ (or is 0 in the case $k=0$).
what remains is to apply the isomorphism of vector spaces $D({\mathbb U})_q
\xrightarrow[\widetilde{\Theta}]{}{\rm End}({\cal H})_f \simeq \bigoplus
\limits_{i \ge 0,\;j \ge 0}{\cal H}_i \otimes{\cal H}_j^*\simeq \bigoplus
\limits_{i \ge 0,\;j \ge 0}{\cal H}_i \otimes{\cal H}_j$. \hfill
$\blacksquare$

\bigskip

\section{Differential forms with finite coefficients}

 The results of this section are not used in producing an explicit formula
for the Bergman kernel. However, it makes an independent interest and
provides an essential addition to the results of \cite{SSV}. In \cite{SV}
a covariant algebra $\Omega({\rm Mat}_{mn})_q$ of differential forms with
polynomial coefficients was considered (it was denoted there by
$\Omega(\mathfrak{g}_{-1})_q$).

 One can find in \cite[section 4]{SSV} a complete list of relations between
the generators $z_a^\alpha$, $(z_a^\alpha)^*$, $dz_a^\alpha$,
$d(z_a^\alpha)^*$, $a=1,\dots,n$, $\alpha=1,\dots,m$, of $\Omega({\rm
Mat}_{mn})_q$. Consider the subalgebras $\bigwedge_{mn}\subset \Omega({\rm
Mat}_{mn})_q$, $\overline{\bigwedge}_{mn}\subset \Omega({\rm Mat}_{mn})_q$,
generated by $\{dz_a^\alpha \}$, $\{d(z_a^\alpha)^*\}$ respectively. They
are q-analogues of algebras of differential forms with constant
coefficients, and $\dim \bigwedge_{mn}=\dim
\overline{\bigwedge}_{mn}=2^{mn}$. There is a decomposition
$$\Omega({\rm
Mat}_{mn})_q={\bigwedge}_{mn}\otimes{\rm Pol}({\rm Mat}_{mn})_q \otimes
\overline{\bigwedge}_{mn}.$$

 We are interested in considering the space $\Omega({\mathbb U})_q
\overset{\rm def}{=}\bigwedge_{mn}\otimes D({\mathbb U})_q \otimes
\overline{\bigwedge}_{mn}$ of differential forms with finite coefficients
and the space $\Omega_q \overset{\rm def}{=}\bigwedge_{mn}\otimes {\rm
Fun}({\mathbb U})_q \otimes \overline{\bigwedge}_{mn}=\Omega({\rm
Mat}_{mn})_q+\Omega({\mathbb U})_q$. We are going to equip $\Omega_q$ with a
structure of covariant differential algebra and to describe it in terms of
generators and relations.

 Use the above topology in $D({\mathbb U})_q'$ to introduce a topology in
the vector space $\Omega'({\mathbb U})_q=\bigwedge_{mn}\otimes D({\mathbb
U})_q'\otimes \overline{\bigwedge}_{mn}=\bigwedge({\rm Mat_{mn}})_q
\bigotimes_{{\mathbb C}[{\rm Mat}_{mn}]_q}D({\mathbb
U})_q'\bigotimes_{{\mathbb C}[\overline{\rm Mat}_{mn}]_q}
\bigwedge(\overline{\rm Mat}_{mn})_q$ of differential forms whose
coefficients are distributions. The differential $d$ and the structure of
covariant $\Omega({\rm Mat}_{mn})_q$-bimodule are transferred by a
continuity from $\Omega({\rm Mat}_{mn})_q$ onto $\Omega'({\mathbb U})_q$.
(In fact, all the commutation relations involving the differentials
$dz_a^\alpha$, $d(z_a^\alpha)^*$ are purely quadratic). Using proposition
\ref{fc}, it is easy to distinguish the differential forms with finite
coefficients from $\Omega'({\mathbb U})_q$, that is, to prove that
$$\Omega({\mathbb U})_q=\{\omega \in \Omega'({\mathbb U})_q|\;\exists
M:\:{\mathbb C}[\overline{\rm Mat}_{mn}]_{q,-M}\cdot \omega=\omega
\cdot{\mathbb C}[{\rm Mat}_{mn}]_{q,M}=0 \}.$$ This allows one to extend by
a continuity the structure of covariant differential algebra from
$\Omega({\rm Mat}_{mn})_q$ onto $\Omega_q=\Omega({\rm
Mat}_{mn})_q+\Omega({\mathbb U})_q$.

 A complete list of commutation relations between the generators
$z_a^\alpha$, $(z_a^\alpha)^*$, $f_0$, $dz_a^\alpha$, $d(z_a^\alpha)^*$,
$a=1,\dots,n$; $\alpha=1,\dots,m$, of $\Omega_q$, includes (\ref{f0dz}) and
the relations from \cite{SSV}. (\ref{df0}) describes the action of the
differential $d$ onto $f_0$.

\medskip

\begin{proposition} For all $a=1,\dots,n$; $\alpha=1,\dots,m$, one has
\begin{equation}\label{f0dz}f_0d(z_a^\alpha)^*=d(z_a^\alpha)^*f_0,\qquad
f_0dz_a^\alpha=dz_a^\alpha f_0.
\end{equation}
\end{proposition}

\smallskip

 {\bf Proof.} It suffices to consider the first relation. It follows from
the invertibility of R-matrices involved into the commutation relations
between $d(z_a^\alpha)^*$ and $z_a^\alpha$, $(z_a^\alpha)^*$, $a=1,\dots,n$;
$\alpha=1,\dots,m$, that
$$f_0d(z_a^\alpha)^*=\sum_{b=1}^n \sum_{\beta=1}^m d(z_b^\beta)^*\psi_{\beta
a}^{b \alpha},\qquad \psi_{\beta a}^{b \alpha}\in D'({\mathbb U})_q.$$

 Prove that $\psi_{\beta a}^{b \alpha}\cdot z_c^\gamma=0$ for all
$c=1,\dots,n$; $\gamma=1,\dots,m$. In fact, $f_0{\mathbb C}[{\rm
Mat}_{mn}]_{q,1}=0$. Hence,
$$0=f_0d(z_a^\alpha)^*z_c^\gamma=\sum_{b=1}^n \sum_{\beta=1}^m
d(z_b^\beta)^*(\psi_{\beta a}^{b \alpha}\cdot z_c^\gamma)$$
(the first equality is due to the homogeneity of the commutation relations
between $d(z_b^\beta)^*$ and $z_c^\gamma$, $a,c=1,\dots,n$;
$\alpha,\gamma=1,\dots,m$). Thus, it follows from the latter equality and
the definition of $\Omega'({\mathbb U})_q$ that $\psi_{\beta a}^{b
\alpha}\cdot z_c^\gamma=0$, $c=1,\dots,n$; $\gamma=1,\dots,m$. \hfill
$\blacksquare$

\medskip

\begin{lemma}\label{ynlr} If $\psi \in D'({\mathbb U})_q$ is such that $\psi
\cdot z_c^\gamma=0$, $c=1,\dots,n$; $\gamma=1,\dots,m$, then $\psi \in
\mathop{\boldsymbol \times}\limits_{j=0}^\infty{\mathbb C}[{\rm
Mat}_{mn}]_{q,j} \cdot f_0 \subset D'({\mathbb U})_q$.
\end{lemma}

\smallskip

 {\bf Proof.} Show first that for any $\psi \in D'({\mathbb U})_q$ one has
$\psi y^N \underset{N \to \infty}{\to}\psi f_0$ in the topology of
$D'({\mathbb U})_q$. In fact, it suffices to demonstrate that $l_
\varphi(\psi y^N)\underset{N \to \infty}{\to}l_ \varphi(\psi f_0)$ for any
$\varphi \in D'({\mathbb U})_q$. Apply the decomposition $D({\mathbb
U})_q=\bigoplus \limits_{\genfrac{}{}{0pt}{1}{j \ge 0}{k \le 0}}D({\mathbb
U})_{q,j,k}$, with $D({\mathbb U})_{q,j,k}={\mathbb C}[{\rm
Mat}_{mn}]_{q,j}\cdot f_0 \cdot{\mathbb C}[\overline{\rm Mat}_{mn}]_{q,k}$.
One has: $\varphi=\sum \limits_{j,k}\varphi_{jk}$, with $\varphi_{jk}\in
D({\mathbb U})_{q,j,k}$, $$l_ \varphi(\psi y^N)=\int \limits_{{\mathbb
U}_q}\psi y^N \varphi d \nu=\sum_{j,k}\int \limits_{{\mathbb U}_q}\psi y^N
\varphi_{jk}d \nu=\sum_{j,k}q^{2Nj}\int \limits_{{\mathbb U}_q}\psi
\varphi_{jk}d \nu.$$ On the other hand, $$\lim_{N \to
\infty}\sum_{j,k}q^{2Nj}\int \limits_{{\mathbb U}_q}\psi \varphi_{jk}d
\nu=\sum_k \int \limits_{{\mathbb U}_q}\psi \varphi_{0k}d \nu=\sum_{j,k}\int
\limits_{{\mathbb U}_q}\psi f_0 \varphi_{jk}d \nu=\int \limits_{{\mathbb
U}_q}\psi f_0 \varphi d \nu=l_ \varphi(\psi f_0).$$

 Turn back to the proof of lemma \ref{ynlr}. If $\psi \in D'({\mathbb U})_q$
and $\psi \cdot z_c^\gamma=0$ for all $c=1,\dots,n$; $\gamma=1,\dots,m$,
then it follows from theorem \ref{yexp} that $\psi y=\psi$, and hence
$\psi=\psi y^N \underset{N \to \infty}{\to}\psi f_0$, so the statement of
the lemma is proved. \hfill $\blacksquare$

\medskip

 We have demonstrated that $\psi \in \mathop{\boldsymbol
\times}\limits_{j=0}^\infty{\mathbb C}[{\rm Mat}_{mn}]_{q,j} \cdot f_0$
($a,b=1,\dots,n$; $\alpha,\beta=1,\dots,m$). Prove that $\psi_{\beta a}^{b
\alpha}$ differs from $f_0$ only by a constant multiple. All $U_q
\mathfrak{sl}_N$-modules considered in \cite{SSV} were equipped with a
gradation determined by the element $H_0$ of the Cartan subalgebra:
$$H_0 \overset{\rm def}{=}\frac{2}{m+n}\left(m \sum_{j=1}^{n-1}jH_j+n
\sum_{j=1}^{m-1}jH_{N-j}+mnH_n \right).$$
An application of the operator $q^{H_0}$ to $f_0d(z_a^\alpha)^*$ and
$d(z_b^\beta)^*\psi_{\beta a}^{b \alpha}$ yields $q^{H_0}(\psi_{\beta a}^{b
\alpha})=\psi_{\beta a}^{b \alpha}$. What remains is to remind that $\{f \in
\mathop{\boldsymbol \times}\limits_{j=0}^\infty{\mathbb C}[{\rm
Mat}_{mn}]_{q,j} \cdot f_0|\;q^{H_0}f=f \}={\mathbb C}f_0$.

 We have proved that the linear span of $\{f_0 \cdot
d(z_a^\alpha)^*\}_{a=1,\dots,n;\;\alpha=1,\dots,m}$ coincides with the
linear span of $\{d(z_a^\alpha)^*f_0 \}_{b=1,\dots,n;\;\beta=1,\dots,m}$.
This vector space is a simple $U_q \mathfrak{sl}_n \otimes U_q
\mathfrak{sl}_m$-module, and the linear map $f_0 \cdot
d(z_a^\alpha)^*\mapsto d(z_a^\alpha)^*f_0$ is the unique up to a constant
multiple endomorphism of this $U_q \mathfrak{sl}_n \otimes U_q
\mathfrak{sl}_m$-module. Hence, $f_0 \cdot d(z_a^\alpha)^*=C \cdot
d(z_a^\alpha)^*f_0$ for some $C \in{\mathbb C}$. On the other hand,
$f_0^*=f_0$, and so $C \cdot d(z_a^\alpha)^*f_0=f_0 \cdot
d(z_a^\alpha)^*=f_0^2 d(z_a^\alpha)^*=C^2 \cdot d(z_a^\alpha)^*f_0^2=C^2
\cdot d(z_a^\alpha)^*f_0$. Therefore, $(1-C)C \cdot d(z_a^\alpha)^*f_0=0$.
On the other hand, $f_0 \cdot d(z_a^\alpha)^*=C \cdot d(z_a^\alpha)^*f_0
\ne 0$, and thus we get $C=1$. \hfill $\blacksquare$

\medskip

\begin{proposition}
\begin{equation}\label{df0}df_0=-\frac{1}{1-q^2}\sum_{a=1}^n
\sum_{\alpha=1}^m(dz_a^\alpha f_0(z_a^\alpha)^*+z_a^\alpha
f_0(dz_a^\alpha)^*).
\end{equation}
\end{proposition}

\smallskip

 {\bf Proof.} It follows from $f_0^*=f_0$ that
$df_0=\overline{\partial}f_0+\partial
f_0=\overline{\partial}f_0+(\overline{\partial}f_0)^*=
\overline{\partial}f_0+(\overline{\partial}f_0)^*$. Hence, it suffices to
prove the relation
$$\overline{\partial}f_0=\omega_0,\qquad
\omega_0=-\frac{1}{1-q^2}\sum_{a=1}^n\sum_{\alpha=1}^mz_a^\alpha
f_0d(z_a^\alpha)^*.$$

 Let $\Omega(\mathbb{U})_q^{(0,1)}$ be the space of (0,1)-forms with finite
coefficients in the quantum ball. Remind that all the $U_q
\mathfrak{sl}_N$-modules in our consideration are equipped with a
$\mathbb{Z}$-grading defined by the element $H_0$ \cite{SSV}. We are about
to prove that 1-forms $\overline{\partial} f_0$ and $\omega_0$ are solutions
of the following system of equations:
\begin{equation}\label{of0cr}(1-y)\omega=\overline{\partial}y \cdot f_0,\qquad
H_0\omega=0
\end{equation}
and to elaborate the uniqueness of a solution of this system in the space
$\Omega(\mathbb{U})_q^{(0,1)}$.

 We start with proving the uniqueness.

\medskip

\begin{lemma} If $\omega \in \Omega(\mathbb{U})_{q}^{(0,1)}$ and
$(1-y)\omega=H_0\omega=0$, then $\omega=0$.
\end{lemma}

\smallskip

 {\bf Proof.} Apply the decomposition $D(\mathbb{U})_q=\bigoplus \limits_{j
\ge 0,k \le 0}D(\mathbb{U})_{q,j,k}$, with $D(\mathbb{U})_{q,j,k}=
\mathbb{C}[\mathrm{Mat}_{mn}]_{q,j}f_0
\mathbb{C}[\overline{\mathrm{Mat}}_{mn}]_{q,k}$. One has $\omega=\sum
\limits_{j=0}^\infty \sum \limits_{k=0}^\infty \omega_{jk}$, where
$\omega_{jk}=\sum \limits_{a=0}^n \sum \limits_{\alpha=0}^mf_{\alpha
jk}^ad(z_a^\alpha)^*$, \ $f_{\alpha jk}\in D(\mathbb{U})_{q,j,-k}$. The
statement of the lemma now follows from the relations $H_0
\omega_{jk}=2(j-k-1)\omega_{jk}$ and
$(1-y)\omega_{jk}=(1-q^{2j})\omega_{jk}$. (The latter relation is
deducible from $yf_0=f_0$ and $yz_a^\alpha=q^2z_a^\alpha y$; $a=1,\ldots,n$;
$\alpha=1,\ldots,m$). \hfill $\blacksquare$

\medskip

 Turn back to the proof of proposition \ref{df0}. What remains is to verify
that the 1-forms $\overline{\partial} f_0$ and $\omega_0$ are solutions of
the equation system \ref{of0cr}.

 The relations $H_0(\overline{\partial} f_0)=0$, $H_0\omega_0=0$ follow from
$H_0f_0=0$, $H_0z_a^\alpha=2z_a^\alpha$,
$H_0(dz_a^\alpha)^*=-2(dz_a^\alpha)^*$, $\alpha=1,\dots,m$, $a=1,\dots,n$,
together with the fact that $\overline{\partial}$ is a morphism of $U_q
\mathfrak{sl}_N$-modules.

 To prove the relation $(1-y)\overline{\partial}f_0=\overline{\partial}y
\cdot f_0$ it is sufficient to apply $\overline{\partial}$ to both sides
of the equality $y \cdot f_0=f_0$.

 The relation $(1-y) \cdot \omega_{0}=\overline{\partial}y \cdot
f_0$ follows from \begin{equation}\label{dyf0}\overline{\partial}y \cdot
f_0=-\sum_{a=1}^n\sum_{\alpha=1}^mz_a^\alpha f_0d(z_a^\alpha)^*.
\end{equation}
On the other hand, (\ref{dyf0}) is deducible from the explicit formula for
the element $y$ obtained in section 2:
$$y=\sum_{k=0}^\infty(-1)^ky_k,\qquad y_0=1,
\qquad y_1=\sum_{a=1}^n\sum_{\alpha=1}^mz_a^\alpha(z_a^\alpha)^*,$$
$$y_k=\sum_{\{J'|\mathrm{card}(J')=k\}}\sum_{\{J''|\mathrm{card}(J'')=k\}}
{z^{\wedge k}}_{J''}^{J'}\left({z^{\wedge k}}_{J''}^{J'}\right)^*.$$ In
fact, $d(z_a^\alpha)^*f_0=f_0d(z_a^\alpha)^*$; $(z_a^\alpha)^*f_0=0$ for all
$a=1,\ldots,n$; $\alpha=1,\ldots,m$, and hence $\overline{\partial}y_k \cdot
f_0=0$ for all $k \ge 2$. \hfill $\blacksquare$

\bigskip

\section{The linear map \boldmath $\stackrel{\circ}{P}_ \lambda$}

 We assume in what follows $\lambda>N-1$.

\medskip

\begin{proposition} The linear subspace $D(\mathbb{U})_q$ is dense in the
Hilbert space $L^2(d \nu_ \lambda)_q$, and the embedding $D(\mathbb{U})_q
\hookrightarrow D(\mathbb{U})_q'$ extends by a continuity up to an embedding
$L^2(d \nu_ \lambda)_q \hookrightarrow D(\mathbb{U})_q'$.
\end{proposition}

\smallskip

 {\bf Proof.} It suffices to apply the relations (\ref{ii}), (\ref{Gam}),
and the isomorphisms
$$D(\mathbb{U})_q \stackrel{\displaystyle
\sim}{\to}\mathrm{End}(\mathcal{H})_f,\qquad
D(\mathbb{U})_q'\overset{\displaystyle
\sim}{\to}\overline{\mathrm{End}}(\mathcal{H})\qquad \text{(see section
4)}\eqno \blacksquare$$

\medskip

 We identify in the sequel the Hilbert space $L^2(d \nu_ \lambda)_q$ and its
image under the embedding into the space of distributions
$D(\mathbb{U})_q'$.

 Consider the orthogonal projection $P_ \lambda$ in the Hilbert space $L^2(d
\nu_ \lambda)_q$ onto the Hardy-Bergman subspace $L_a^2(d \nu_ \lambda)_q$
introduced in section 3. A principal subject of the research in the
remainder of this work will be the linear map $\stackrel{\circ}{P}_
\lambda:D(\mathbb{U})_q \to D(\mathbb{U})_q'$ given by a restriction of $P_
\lambda$ onto the dense in $L^2(d \nu_ \lambda)_q$ linear subspace
$D(\mathbb{U})_q$. \footnote{It will be proved in the sequel that $P_
\lambda \mathrm{Fun}(\mathbb{U})_q \subset \mathbb{C}[\mathrm{Mat}_{mn}]_q$}

 This section presents a construction of such a representation $\pi_
\lambda$ of $U_q \mathfrak{sl}_N$ in $D(\mathbb{U})_q'$ that $\pi_
\lambda(a)D(\mathbb{U})_q \subset D(\mathbb{U})_q$ for all $a \in U_q
\mathfrak{sl}_N$, and
\begin{equation}\label{plcr}\pi_ \lambda(a)\stackrel{\circ}{P}_
\lambda=\stackrel{\circ}{P}_ \lambda\pi_ \lambda(a),\qquad a \in U_q
\mathfrak{sl}_N.
\end{equation}

 The results of section 7 will imply that $P_ \lambda$ is uniquely
determined by (\ref{plcr}) and its value on $f_0$.

\medskip

\begin{proposition}\label{pl} There exists a unique representation $\pi_
\lambda$ of $U_q \mathfrak{sl}_N$ in $D(\mathbb{U})_q'$ such that for all $f
\in D(\mathbb{U})_q'$
\begin{equation}\label{pl1}\pi_ \lambda(E_j):f \mapsto
\begin{cases} E_jf, & j \ne n, \\ E_nf-q^{1/2}\dfrac{1-q^{2
\lambda}}{1-q^2}(K_nf)z_n^m, & j=n.
\end{cases}
\end{equation}
\begin{equation}\pi_ \lambda(F_j):f \mapsto
\begin{cases} F_jf, & j \ne n, \\ q^{-\lambda} F_nf, & j=n.
\end{cases}
\end{equation}
\begin{equation}\label{pl3}\pi_ \lambda(K_j^{\pm 1}):f \mapsto
\begin{cases} K_j^{\pm 1}f, & j \ne n, \\ q^{\pm \lambda}K_n^{\pm 1}f, & j=n.
\end{cases}
\end{equation}
\end{proposition}

\smallskip

 {\bf Proof.} The uniqueness of $\pi_ \lambda$ is obvious. While proving the
existence of this representation, it suffices to replace the topological
vector space $D(\mathbb{U})_q'$ by its dense subspace
$\mathrm{Pol}(\mathrm{Mat}_{mn})_q$, and to consider the special case
$\lambda \in \{N,N+1,N+2,\ldots \}$. (In fact, the problem is to prove
equalities in which both sides are in $\mathbb{C}[q^\lambda,q^{-\lambda}]$.
So, what remains is to observe that two polynomials which coincide on the
set $\{q^N,q^{N+1},q^{N+2}\ldots \}$ are identically the same.) Consider the
$U_q \mathfrak{sl}_N$-module $\mathrm{Pol}(\widetilde{X})_{q,x}$ and the
associated representation $\pi$ of $U_q \mathfrak{sl}_N$. Remind the
notation $t=t_{\{1,2,\dots,m \}\{n+1,n+2,\dots,N \}}^{\wedge m}$. The
existence of $\pi_ \lambda$ follows from the following

\medskip

\begin{lemma}\label{pl0} Let $\lambda \in \{N,N+1,N+2,\ldots \}$ and
${\EuScript I}_{\lambda}$ be the linear map ${\EuScript
I}_{\lambda}:\mathrm{Pol}(\mathrm{Mat}_{mn})_q \to
\mathrm{Pol}(\widetilde{X})_{q,x};\qquad{\EuScript I}_{\lambda}:f
\mapsto ({\EuScript
I}f)t^{-\lambda}$. Then for all $j=1,\ldots,N-1$ one has $\pi_
\lambda(E_j)={\EuScript I}_{\lambda}^{-1}\pi(E_j){\EuScript I}_{\lambda}$,
\ $\pi_
\lambda(F_j)={\EuScript I}_{\lambda}^{-1}\pi(F_j){\EuScript I}_{\lambda}$,
 \ $\pi_
\lambda(K_j^{\pm 1})={\EuScript I}_{\lambda}^{-1}\pi(K_j^{\pm 1}){\EuScript I}
_{\lambda}$.
\end{lemma}

\smallskip

 {\bf Proof.} It follows from the covariance of the algebra
$D(\mathbb{U})_q'$ that
$$E_j(ft^{-\lambda})=(E_jf)(t^{-\lambda})+(K_jf)(E_j(t^{-\lambda})),$$
$$F_j(ft^{-\lambda})=(F_jf)(K_j^{-1}(t^{-\lambda}))+f
\cdot(F_j(t^{-\lambda})),$$
$$K_j^{\pm 1}(ft^{-\lambda})=(K_j^{\pm 1}f)(K_j^{\pm 1}(t^{-\lambda})),$$
so it suffices to prove the relations
$$F_j(t^{-\mu})=0,\qquad K_j^{\pm 1}(t^{-\mu})=
\begin{cases} q^{\pm \mu}t^{-\mu}, & j=n \\ t^{-\mu}, & j \ne n.
\end{cases}$$
$$E_j(t^{-\mu})=
\begin{cases}-q^{1/2}\dfrac{1-q^{2 \mu}}{1-q^2}\left(t^{-1}\cdot
t_{\{1,2,\dots,m \}\{n,n+2,\dots,N \}}^{\wedge m}\right)t^{-\mu}, & j=n,\\
0, & j \ne n.
\end{cases}$$
for all $\mu \in-\mathbb{Z}_+$. These are easily deducible via an
application of the covariance of $\mathrm{Pol}(\widetilde{X})_{q,x}$ and the
relations
$$F_jt=0,\qquad K_j^{\pm 1}t=\begin{cases} q^{\mp 1}t,& j=n\\ t, & j
\ne n.\end{cases},$$
$$E_jt=\begin{cases}q^{-1/2}\cdot t_{\{1,2,\dots,m \}\{n,n+2,\dots,N
\}}^{\wedge m}, & j=n\\ 0, & j \ne n.\end{cases}.$$

 (For example,
$$E_n(t^k)=\sum_{j=0}^{k-1}(K_nt)^j(E_nt)t^{k-j-1}=
q^{-1/2}\sum_{j=0}^{k-1}(q^{-1}t)^j t_{\{1,2,\dots,m \}\{n,n+2,\dots,N
\}}^{\wedge m}t^{k-j-1}=$$
$$=q^{-3/2}(\sum_{j=0}^{k-1}q^{-2j})\left(t^{-1}t_{\{1,2,\dots,m
\}\{n,n+2,\dots,N \}}^{\wedge m}\right)t^k.)\eqno \blacksquare$$

 Remind the notion of an invariant scalar product in a representation space
of a Hopf $*$-algebra $A$. Let $S$ be the antipode and $\varepsilon$ the
counit of this Hopf algebra.

 Consider an $A$-module $V$. The antimodule $\overline{V}$ is defined to be
$V$ as an Abelian group, while the multiplication by complex numbers and
$A$-action in $\overline{V}$ are given by
$$(\lambda,v)\mapsto \overline{\lambda}v,\qquad (a,v)\mapsto(S(a))^*v,\qquad
\lambda \in \mathbb{C},\;a\in A,\;v\in \overline{V}.$$

 Let $V_1$, $V_2$ be two $A$-modules. A sesquilinear form $V_1 \times V_2
\to \mathbb{C}$ is called invariant if the associated linear functional
$\eta:\overline{V_2}\otimes V_1 \to \mathbb{C}$ is a morphism of
$A$-modules:
$$\eta(av)=\varepsilon(a)\eta(v),\qquad a\in A,\;v \in
\overline{V_2}\otimes V_1.$$

 Consider the representation $\pi_ \lambda$ and its subrepresentation in
$D(\mathbb{U})_q$. Let $D(\mathbb{U})_{q,\lambda}$ be the associated $U_q
\mathfrak{su}_{n,m}$-module. (The $U_q \mathfrak{su}_{n,m}$-modules
$\mathbb{C}[\overline{\mathrm{Mat}}_{mn}]_{q,\lambda}$,
$\mathrm{Pol}(\overline{\mathrm{Mat}}_{mn})_{q,\lambda}$, and
$\mathrm{Fun}(\mathbb{U})_{q,\lambda}=
\mathrm{Pol}(\overline{\mathrm{Mat}}_{mn})_{q,\lambda}+
D(\mathbb{U})_{q,\lambda}$ are defined in a similar way.)

\medskip

\begin{proposition} For all $\lambda>N-1$, the scalar product
$D(\mathbb{U})_{q,\lambda}\times D(\mathbb{U})_{q,\lambda}\to \mathbb{C}$,
$f_1 \times f_2 \mapsto \displaystyle\int \limits_{\mathbb{U}_q}f_2^*\cdot
f_1y^\lambda d \nu$ is positive and $U_q \mathfrak{su}_{n,m}$-invariant.
\end{proposition}

\smallskip

 {\bf Proof.} The positivity was demonstrated before (in section 3). The
same argument as in the proof of proposition \ref{pl} allows one to reduce
matters to the special case $\lambda \in \{N,N+1,N+2,\ldots \}$.

 In this special case one has a well defined operator $\EuScript{I}_
\lambda:D(\mathbb{U})_{q,\lambda}\to D(\widetilde{X})_q$, $\EuScript{I}_
\lambda:f \mapsto (\EuScript{I}f)t^{-\lambda}$. It follows from
lemma \ref{pl0} that this linear map is a morphism of $U_q
\mathfrak{su}_{n,m}$-modules.

 One can find in \cite{SSV} a construction of invariant integral on a
quantum principal homogeneous space. It is easy to deduce from that
construction that $\displaystyle \int \limits_{\widetilde{X}_q}ft^*d
\nu=\displaystyle \int \limits_{\widetilde{X}_q}t^*fd \nu$ for all $f \in
D(\widetilde{X})_q$ since the operator $\Pi(tt^*)$ introduced in the paper
alluded above commutes with $\Gamma(e^{h \check{\rho}})$ involved in
(\ref{ii}). Hence
$$\int \limits_{\mathbb{U}_q}f_2^*f_1y^\lambda d \nu=\int
\limits_{\widetilde{X}_q}(\EuScript{I}f_2)^*(\EuScript{I}f_1)
(tt^*)^{-\lambda}d \nu=\int
\limits_{\widetilde{X}_q}(\EuScript{I}f_2)^*(\EuScript{I}f_1)
t^{-\lambda}(t^*)^{-\lambda}d \nu=$$
$$=\int
\limits_{\widetilde{X}_q}(t^*)^{-\lambda}(\EuScript{I}f_2)^*(\EuScript{I}f_1)
t^{-\lambda}d \nu=\int
\limits_{\widetilde{X}_q}(\EuScript{I}_{\lambda}f_2)^*(\EuScript{I}_{\lambda}
f_1)d
\nu.$$ Thus, the invariance of scalar product in $D(\mathbb{U})_{q,\lambda}$
follows from the invariance of the scalar product $D(\widetilde{X})_q \times
D(\widetilde{X})_q \to \mathbb{C}$, $f_1 \times f_2 \mapsto \int
\limits_{\widetilde{X}_q}f_2^*f_1d \nu$, while the latter statement follows
from the invariance of the integral on the quantum principal homogeneous
space.

\medskip

\begin{corollary}\label{spi} For all $\lambda>N-1$, the scalar product
$\mathrm{Fun}(\mathbb{U})_{q,\lambda}\times
\mathrm{Fun}(\mathbb{U})_{q,\lambda}\to \mathbb{C}$, $f_1 \times f_2
\mapsto(f_1,f_2)_\lambda \stackrel{\rm def}{=}C(\lambda)\displaystyle \int
\limits_{\mathbb{U}_q}f_2^*f_1d \nu_ \lambda$ is $U_q
\mathfrak{su}_{n,m}$-invariant (the constant $C(\lambda)$ is determined by
(\ref{Cl})).
\end{corollary}

\smallskip

 {\bf Proof.} Let $j \in \mathbb{Z}_+$ and $\chi_j \in D(\mathbb{U})_q$ be
such a finite function that the operator $\Theta(\chi_j)$ in
$\mathcal{H}=\bigoplus_{k=0}^{\infty} \mathcal{H}_k$ is the projection
onto the space $\bigoplus_{k=0}^j \mathcal{H}_k$ along the subspace
$\bigoplus_{k=j+1}^\infty \mathcal{H}_k$.  (The existence and uniqueness of
such $\chi$ is due to the isomorphism $\Theta:D(\mathbb{U})_q \to
\mathrm{End}(\mathcal{H})_f$). The invariance of the scalar product
$(f_1,f_2)_ \lambda$ in $D(\mathbb{U})_{q,\lambda}$ implies the invariance
of the associated scalar product in $\mathrm{Fun}(\mathbb{U})_{q,\lambda}$
since for all $f_1,f_2 \in \mathrm{Fun}(\mathbb{U})_{q,\lambda}$, $a_1,a_2
\in U_q \mathfrak{su}_{n,m}$ $$(\pi_ \lambda(a_1)f_1,\pi_ \lambda(a_2)f_2)_
\lambda=\lim_{m_1,m_2 \to \infty}(\pi_ \lambda(a_1)(\chi_{m_1}f_1
\chi_{m_2}),\pi_ \lambda(a_2)(\chi_{m_2}f_2 \chi_{m_1}))_ \lambda.$$ While
proving the latter equality, one should use the description of $U_q
\mathfrak{sl}_N$-action in $D(\mathbb{U})_q$ from \cite[sections 7, 8]{SSV}.
\hfill $\blacksquare$

\medskip

\begin{proposition}\label{Pl} $P_ \lambda\,\mathrm{Fun}(\mathbb{U})_{q,
\lambda}\subset
\mathbb{C}[\mathrm{Mat}_{mn}]_{q,\lambda}$, and the associated linear map
$P_ \lambda:\mathrm{Fun}(\mathbb{U})_{q,\lambda}\to
\mathbb{C}[\mathrm{Mat}_{mn}]_{q,\lambda}$ is a morphism of $U_q
\mathfrak{su}_{n,m}$-modules for all $\lambda>N-1$.
\end{proposition}

\smallskip

 {\bf Proof.} Each $f \in \mathrm{Fun}(\mathbb{U})_{q,\lambda}$ is
orthogonal to all but finitely many of homogeneous components of the graded
vector space $\mathbb{C}[\mathrm{Mat}_{mn}]_{q,\lambda}$. Hence $P_ \lambda
f \in \mathbb{C}[\mathrm{Mat}_{mn}]_{q,\lambda}$. Since $P_ \lambda^2=P_
\lambda$, it suffices to prove that
$\mathbb{C}[\mathrm{Mat}_{mn}]_{q,\lambda}$ and its orthogonal complement in
$\mathrm{Fun}(\mathbb{U})_{q,\lambda}$ are $U_q
\mathfrak{su}_{n,m}$-submodules of the $U_q \mathfrak{su}_{n,m}$-module
$\mathrm{Fun}(\mathbb{U})_{q,\lambda}$. For the first subspace this follows
from the definition of $\pi_ \lambda$, and for the second one this fact is
due to corollary \ref{spi} (by the invariance of the scalar product we have:
$$(\pi_\lambda(a)f_1,f_2)_ \lambda=(f_1,\pi_\lambda(a^*)f_2)_ \lambda,\qquad
a \in U_q \mathfrak{su}_{n,m},\;f_1,f_2 \in
\mathrm{Fun}(\mathbb{U})_{q,\lambda}$$ (cf. \cite{SSV3}). \hfill
$\blacksquare$

\medskip

\begin{corollary} $P_ \lambda\,D(\mathbb{U})_{q,\lambda}\subset
\mathbb{C}[\mathrm{Mat}_{mn}]_{q,\lambda}$, and the operator
$\stackrel{\circ}{P}_ \lambda:D(\mathbb{U})_{q,\lambda}\to
\mathbb{C}[\mathrm{Mat}_{mn}]_{q,\lambda}$ is a morphism of $U_q
\mathfrak{su}_{n,m}$-modules.
\end{corollary}

\bigskip

\section{The element $f_0$}

Consider the subalgebras $U_q \mathfrak{N}_ \pm \subset U_q \mathfrak{sl}_N$
generated by $\{E_j \}_{j=1,\ldots,N-1}$ and $\{F_j \}_{j=1,\ldots,N-1}$,
respectively.

\medskip

\begin{lemma}\label{ff} $f_0$ generates the $U_q \mathfrak{N}_+$-module
$\mathbb{C}[\mathrm{Mat}_{mn}]_qf_0$ and the $U_q \mathfrak{N}_-$-module
$f_0\mathbb{C}[\overline{\mathrm{Mat}}_{mn}]_q$.
\end{lemma}

\smallskip

 {\bf Proof.} It suffices to prove the first statement. One can find in
\cite{SV, SSV} a description of the generalized Verma module $V_+(0)$ over
$U_q \mathfrak{sl}_N$ with a single generator $v_+(0)\in V_+(0)$. It is
sufficient to demonstrate that the map $v_+(0)\mapsto f_0$ admits an
extension up to an isomorphism of the {\sl graded} $U_q
\mathfrak{N}_+$-modules $V_+(0)\overset{\displaystyle
\sim}{\to}\mathbb{C}[\mathrm{Mat}_{mn}]_qf_0$. On the other hand,
$\mathbb{C}[\mathrm{Mat}_{mn}]_qf_0$ is a dual {\sl graded} $U_q
\mathfrak{N}_+$-module with respect to
$\mathbb{C}[\overline{\mathrm{Mat}}_{mn}]_q$, due to the invariance and
nondegeneracy of the bilinear form $\mathbb{C}[\mathrm{Mat}_{mn}]_qf_0
\times\mathbb{C}[\overline{\mathrm{Mat}}_{mn}]_q \to\mathbb{C}$; $f_1 \times
f_2\mapsto \int \limits_{\mathbb{U}_q}f_1f_2d \nu$ (nondegeneracy follows
from \cite[lemma 8.4]{SSV}). Furthermore, it was shown in \cite{SV, SSV}
that $V_+(0)\simeq(\mathbb{C}[\overline{\mathrm{Mat}}_{mn}]_q)^*$ in the
category of $U_q \mathfrak{sl}_N$-modules. What remains is to refer to the
coincidence of the kernels of the linear functionals associated to $f_0$ and
$v_+(0)$ under the above isomorphisms. (These kernels are just $\bigoplus
\limits_{j=1}^\infty \mathbb{C}[\overline{\mathrm{Mat}}_{mn}]_{q,-j}$).
\hfill $\blacksquare$

\medskip

\begin{proposition}\label{f0g} $U_q \mathfrak{sl}_Nf_0=D(\mathbb{U})_q$.
\end{proposition}

\smallskip

 {\bf Proof.} Consider an ordered set $\{j_1,j_2,\ldots j_r \}$, $r \in
\mathbb{Z}_+$, formed by the elements of the set $\{1,2,\ldots,N-1 \}$. It
suffices to prove that $(E_{j_1}E_{j_2}\ldots
E_{j_r}f_0)\mathbb{C}[\overline{\mathrm{Mat}}_{mn}]_q\subset U_q
\mathfrak{sl}_Nf_0$, since the linear span of $E_{j_1}E_{j_2}\ldots
E_{j_r}f_0$ coincides with $\mathbb{C}[\mathrm{Mat}_{mn}]_qf_0$ by a virtue
of lemma \ref{ff}. We proceed by induction in $r$. In the case $r=0$ our
statement follows from lemma \ref{ff}. The induction passage from $r-1$ to
$r$ could be easily done via an application of
\begin{equation}\label{ir}(E_j(f_+f_0))f_-=E_j(f_+f_0f_-)-K_j(f_+f_0)(E_jf_-).
\end{equation}
In fact, set up $j=j_1$, $f_+f_0=E_{j_2}E_{j_3}\ldots E_{j_r}f_0$, $f_-\in
\mathbb{C}[\overline{\mathrm{Mat}}_{mn}]_q$. By the induction hypothesis one
has
$$f_+f_0f_-\in
(E_{j_2}E_{j_3}\ldots
E_{j_r}f_0)\mathbb{C}[\overline{\mathrm{Mat}}_{mn}]_q\subset U_q
\mathfrak{sl}_Nf_0,$$
$$K_{j_1}(f_+f_0)(E_{j_1}f_-)\in
(E_{j_2}E_{j_3}\ldots E_{j_r}f_0)\mathbb{C}[\overline{\mathrm{Mat}}_{mn}]_q
\subset U_q \mathfrak{sl}_Nf_0.$$
 Hence $(E_{j_1}E_{j_2}\ldots E_{j_r}f_0)f_-=(E_{j_1}(f_+f_0))f_-=
E_{j_1}(f_+f_0f_-)-K_{j_1}(f_+f_0)E_{j_1}f_-\subset U_q \mathfrak{sl}_Nf_0$.
\hfill $\blacksquare$

\medskip

The relations (\ref{pl1}) -- (\ref{pl3}) allow one to generalize the
statements of lemma \ref{ff} and proposition \ref{f0g}.

\medskip

\begin{lemma}\label{pl4} $\{\pi_ \lambda(\xi)f_0|\;\xi \in U_q \mathfrak
{N}_+\}=
\mathbb{C}[\mathrm{Mat}_{mn}]_qf_0$,\\ $\{\pi_ \lambda(\xi)f_0|\;\xi \in U_q
\mathfrak{N}_-\}=f_0 \mathbb{C}[\overline{\mathrm{Mat}}_{mn}]_q$.
\end{lemma}

\smallskip

 {\bf Proof.} The first statement follows from lemma \ref{ff} since the
action of the operators $\pi_\lambda(\xi)$, $\xi \in U_q \mathfrak{N}_+$, on
the subspace $\mathbb{C}[\mathrm{Mat}_{mn}]_qf_0$ is independent of
$\lambda$. The second statement reduces to lemma \ref{ff} via replacement of
the generator $F_n \mapsto q^{-\lambda}F_n$. \hfill $\blacksquare$

\medskip

\begin{proposition}\label{f0gdu} $\{\pi_ \lambda(\xi)f_0|\;\xi \in U_q
\mathfrak{sl}_N\}=D(\mathbb{U})_q$
\end{proposition}

\smallskip

 {\bf Proof.} Repeat the proof of proposition \ref{f0g} with the reference
to lemma \ref{ff} being replaced by that to lemma \ref{pl4}. The statement
$(E_{j_1}E_{j_2}\ldots
E_{j_r}f_0)\mathbb{C}[\overline{\mathrm{Mat}}_{mn}]_q\subset \{\pi_
\lambda(\xi)f_0|\;\xi \in U_q \mathfrak{sl}_N\}$ is proved by induction in
$r$ as before. The only difference is that the first term in the right hand
side of (\ref{ir}) should be replaced in the case $j=n$ by $\pi_
\lambda(E_n)(f_+f_0f_-)+q^{1/2}\dfrac{1-q^{2 \lambda}}{1-q^2}
K_n(f_+f_0f_-)z_n^m$. The appearance of the term ${\tt
const}(f_+,\lambda)f_+f_0K_n(f_-)z_n^m$ does not require introducing any
essential changes to the induction process in question since $K_n(f_-)\in
\mathbb{C}[\overline{\mathrm{Mat}}_{mn}]_q$, $f_+f_0K_n(f_-)z_n^m \in f_+f_0
\mathbb{C}[\overline{\mathrm{Mat}}_{mn}]_q$. \hfill $\blacksquare$

\medskip

 {\sc Remark 7.5.} It follows from proposition \ref{f0g} that the invariant
integral $D(\mathbb{U})_q \to \mathbb{C}$ on the quantum matrix ball is
unique up to a constant multiple.

\bigskip

\section{The integral operators \boldmath $\widetilde{K}_l$}

 This section contains a construction of a family of integral operators
$\widetilde{K}_l:D(\mathbb{U})_q \rightarrow D(\mathbb{U})_q'$ which commute
with the operators of the representation $\pi_l$.

 Let $\mathbb{C}[\mathrm{Mat}_{mn}]_q^{\rm op}$ be the graded algebra
derived from $\mathbb{C}[\mathrm{Mat}_{mn}]_q$ via replacing its
multiplication law with the opposite one. The term 'algebra of kernels' will
stand for a completion of the bigraded algebra
$\mathbb{C}[\mathrm{Mat}_{mn}]_q^{\rm op}\otimes
\mathbb{C}[\overline{\mathrm{Mat}}_{mn}]_q$, that is, the algebra of formal
series of the form $K=\sum \limits_{i,j=0}^\infty K^{(i,j)}$, $ K^{(i,j)}\in
\mathbb{C}[\mathrm{Mat}_{mn}]_{q,i}^{\rm
op}\otimes\mathbb{C}[\overline{\mathrm{Mat}}_{mn}]_{q,-j}$, with the
topology of coefficientwise convergence (the topology of direct product).
This algebra is denoted by $\mathbb{C}[[\mathrm{Mat}_{mn}\times
\overline{\mathrm{Mat}}_{mn}]]_q$.

 To begin with, we construct the kernels $K_l$ of integral operators
$\widetilde{K}_l$ in the special case $l \in-\mathbb{N}$. A passage to the
general case is to be performed later on via an 'analytic continuation' with
respect to the parameter $l$ in $\mathbb{C}[[\mathrm{Mat}_{mn}\times
\overline{\mathrm{Mat}}_{mn}]]_q$ (cf. \cite{SSV3}).

 It follows from the definition of the involutions $*$, $\star$ that
$$\left(t_{\{1,2,\ldots,m \}J}^{\wedge m}\right)^*=(-1)^{({\rm
card}\{1,2,\ldots,n \}\cap J)}\left(t_{\{1,2,\ldots,m \}J}^{\wedge
m}\right)^\star,$$
 where $\left(t_{\{1,2,\ldots,m \}J}^{\wedge
m}\right)^\star=(-q)^{l(J)}t_{\{m+1,\ldots,N \}J'}^{\wedge n}$,
$J'=\{1,\ldots,N \}\setminus J$, $l(J)={\rm card}\{(a,b)|\;a>b\quad \& \quad
a \in J \quad \& \quad b \in J' \}$. Apply these relations and the $U_q
\mathfrak{sl}_N$-invariance of the element $\sum \limits_{s \in
S_N}(-q)^{l(s)}t_{1s(1)}t_{2s(2)}\cdots t_{ms(m)}\otimes
t_{m+1\,s(m+1)}t_{m+2\,s(m+2)}\cdots t_{Ns(N)}\in \mathbb{C}[SL_{N}]_q
\otimes \mathbb{C}[SL_{N}]_q$  to obtain

\medskip

\begin{lemma}\label{L} The 'kernel'
\begin{equation}\label{}L=\sum \limits_{{\rm card}(J)=m,\,J\subset
\{1,2,\ldots,N \}}(-1)^{({\rm card}\{1,2,\ldots,n \}\cap
J)}t_{\{1,2,\ldots,m \}J}^{\wedge m}\otimes \left(t_{\{1,2,\ldots,m
\}J}^{\wedge m}\right)^*
\end{equation}
is a $U_q \mathfrak{sl}_N$-invariant of the $U_q \mathfrak{sl}_N$-module
$\mathrm{Pol}(\widetilde{X})_q\otimes \mathrm{Pol}(\widetilde{X})_q$. (That
is, $aL=\varepsilon(a)L$ for all $a \in U_q \mathfrak{sl}_N$.)
\end{lemma}

\medskip

 Consider the algebra $\mathrm{Pol}(\widetilde{X})_{q,x}^{\rm op}$ which is
coming from $\mathrm{Pol}(\widetilde{X})_{q,x}$ via replacing its
multiplication law with the opposite one.

\medskip \stepcounter{theorem}

 {\sc Remark 8.2.} An application of proposition \ref{Izm} allows one to
prove that in the algebra $\mathrm{Pol}(\widetilde{X})_q^{\rm op}\otimes
\mathrm{Pol}(\widetilde{X})_q$
$$ L={\EuScript I}\otimes{\EuScript
I}\left(1+\sum_{k=1}^m(-1)^k\chi_k \right)t \otimes t^*,$$
 with $\chi_k \in \mathbb{C}[\mathrm{Mat}_{mn}]_q^{\rm op}\otimes
\mathbb{C}[\mathrm{Mat}_{mn}]_q \subset \mathbb{C}[[\mathrm{Mat}_{mn}\times
\overline{\mathrm{Mat}}_{mn}]]_q$ being the kernels given by
\begin{equation}\label{hik}\chi_k=\sum_{\genfrac{}{}{0mm}{1}{J'\subset\{1,
\dots,m
\}}{{\rm card}(J')=k}} \sum_{\genfrac{}{}{0mm}{1}{J''\subset\{1,\dots,n
\}}{{\rm card}(J'')=k}}z_{\quad J''}^{\wedge kJ'}\otimes \left(z_{\quad
J''}^{\wedge kJ'}\right)^*.
\end{equation}

 It was shown in \cite[section 2]{V} that a product of any two $U_q
\mathfrak{sl}_N$-invariants of $\mathrm{Pol}(\widetilde{X})_q^{\rm
op}\otimes \mathrm{Pol}(\widetilde{X})_q$ is again a $U_q
\mathfrak{sl}_N$-invariant. Hence the following generalization of lemma
\ref{L}.

\medskip

\begin{lemma} All the powers $L^j$, $j \in \mathbb{N}$, of $L \in
\mathrm{Pol}(\widetilde{X})_q^{\rm op}\otimes
\mathrm{Pol}(\widetilde{X})_q$, are $U_q \mathfrak{sl}_N$-invariants.
\end{lemma}

\medskip

 Let $l \in-\mathbb{N}$. Define the kernel $K_l$ by
\begin{equation}\label{Kl}K_l=\left(1+\sum_{k=1}^m(-q^{2l})^k \chi_k
\right)\cdot \left(1+\sum_{k=1}^m(-q^{2(l+1)})^k \chi_k \right)\cdots
\left(1+\sum_{k=1}^m(-q^{-2})^k \chi_k \right).
\end{equation}

\medskip

\begin{corollary} For all $l \in-\mathbb{N}$ the element $(t \otimes
t^*)^{-l}\cdot{\EuScript I}\otimes{\EuScript I}(K_l)$ of
$\mathrm{Pol}(\widetilde{X})_q^{\rm op}\otimes
\mathrm{Pol}(\widetilde{X})_q$ is equal to $L^{-l}$ and hence is a $U_q
\mathfrak{sl}_N$-invariant.
\end{corollary}

\smallskip

 {\bf Proof.} It suffices to apply remark 8.2 and the commutation relation
$${\EuScript I}\otimes{\EuScript I}(\chi_k)(t \otimes t^*)=q^{-2k}(t \otimes
t^*){\EuScript I}\otimes{\EuScript I}(\chi_k),$$
 which follows from
$${\EuScript I}(z_a^\alpha)t=qt{\EuScript I}(z_a^\alpha),\qquad
a=1,\ldots,n;\;\alpha=1,\ldots,m.\eqno \blacksquare$$

\medskip

 Consider the integral operator $\widehat{K}_l:D(\mathbb{U})_q\to
D(\mathbb{U})_q'$; $\widehat{K}_l:f \mapsto {\rm id}\otimes \nu(K_l(1
\otimes fy^l))$, with $l \in-\mathbb{N}$, and $\nu:D(\mathbb{U})_q \to
\mathbb{C}$ being an invariant integral.

\medskip

\begin{proposition}\label{plkl} For all $l \in-\mathbb{N}$, $a \in U_q
\mathfrak{sl}_N$ one has the equality of operators from $D(\mathbb{U})_q$
into itself
$$\pi_l(a)\widehat{K}_l=\widehat{K}_l \pi_l(a).$$
\end{proposition}

\smallskip

 {\bf Proof.} Consider the integral operator $\widehat{\EuScript{K}}_l$ on
the quantum principal homogeneous space determined by its kernel
$\EuScript{K}_l=(t \otimes t^*)^{-l}{\EuScript I}\otimes{\EuScript I}(K_l)$.
We need also an extension by a continuity of the map ${\EuScript
I}_l:\mathrm{Pol}(\mathrm{Mat}_{mn})_q \to
\mathrm{Pol}(\widetilde{X})_{q,x}$; ${\EuScript I}_l:f \mapsto{\EuScript
I}(f)t^{-l}$ onto the space $\mathrm{Fun}(\mathbb{U})_q=
\mathrm{Pol}(\mathrm{Mat}_{mn})_q+D(\mathbb{U})_q$.

 It suffices to prove the relations
$$\pi(a)\widehat{\EuScript{K}}_lf=\widehat{\EuScript{K}}_l \pi(a)f,\qquad a
\in U_q \mathfrak{sl}_N,\;f \in D(\widetilde{X})_q,$$
$$\pi_{l}(a)={\EuScript I}_l^{-1}\pi(a){\EuScript I}_l,\qquad a \in U_q
\mathfrak{sl}_N,$$
$$\widehat{K}_l={\EuScript I}_l^{-1}\widehat{\EuScript{K}}_l{\EuScript I}_l.$$
The first of those follows from the invariance of the kernel
${\EuScript{K}}_l$ and invariance of the integral involved when constructing
the operator $\widehat{\EuScript{K}}_l$. The second relation is a
consequence of lemma \ref{pl0} and a continuity argument. The latter
equality follows from the fact that for all $\psi \in
\mathrm{Pol}(\mathrm{Mat}_{mn})_q$ and $f \in D(\mathbb{U})_q$ one has
$\displaystyle \int
\limits_{\widetilde{X}_q}(t^*)^{-l}\EuScript{I}(\psi)\EuScript{I}(f)t^{-l}d
\nu=const \int \limits_{\mathbb{U}_q}\psi fd \nu_l$. This relation is proved
as follows:  $$\int
\limits_{\widetilde{X}_q}(t^*)^{-l}\EuScript{I}(\psi)\EuScript{I}(f)t^{-l}d
\nu=\int \limits_{\widetilde{X}_q}(t^*)^{-l}\EuScript{I}(\psi f)t^{-l}d
\nu=\int \limits_{\widetilde{X}_q}\EuScript{I}(\psi f)t^{-l}(t^*)^{-l}d
\nu=\int \limits_{\widetilde{X}_q}\EuScript{I}(\psi f)(tt^*)^{-l}d \nu=$$
$$=\int \limits_{\widetilde{X}_q}\EuScript{I}(\psi fy^l)d \nu=\int
\limits_{\mathbb{U}_q}\psi fy^ld \nu=const \int \limits_{\mathbb{U}_q}\psi fd
\nu_l.$$ (The above argument applies the relations $tt^*=t^*t$ and
$\displaystyle \int \limits_{\widetilde{X}_q}t^*fd \nu=\int
\limits_{\widetilde{X}_q}ft^*d \nu$, $f \in D(\widetilde{X})_q$, together
with theorem \ref{yexp}. The latter relation follows from the special case
$f \in D(\widehat{X})_q$, and hence from even more special case $f=\varphi
t$, $\varphi \in D(X)_q$. Now for $\varphi \in D(X)_q$ one has
$\displaystyle \int \limits_{\widetilde{X}_q}t^* \varphi td \nu=\int
\limits_{\widetilde{X}_q}\varphi tt^*d \nu$, as one can easily deduce from
the explicit formula for invariant integral.) \hfill $\blacksquare$

\medskip

 Now pass from the special case $l \in-\mathbb{N}$ to the general case via
'analytic continuation'.

 We need in the sequel some integral operators whose kernels depend on a
parameter $u$.

 The term 'polynomial kernels' will stand for the formal series $K(u)=\sum
\limits_{i,j=0}^\infty K(u)^{(i,j)}$ whose terms belong to the
$\mathbb{C}[u]$-module $\mathbb{C}[\mathrm{Mat}_{mn}]_{q,i}^{\rm op}\otimes
\mathbb{C}[\mathrm{Mat}_{mn}]_{q,-j}\otimes \mathbb{C}[u]$. The vector space
of all polynomial kernels carries a natural structure of algebra over
$\mathbb{C}[u]$. For any $u_0 \in \mathbb{C}$ one has a well defined
homomorphism $K(u)\mapsto K(u_0)$ of this algebra into
$\mathbb{C}[[\mathrm{Mat}_{mn}\times \overline{\mathrm{Mat}}_{mn}]]_q$.

\medskip

\begin{proposition}\label{pk} There exists a unique polynomial kernel $K(u)$
such that for all $l \in-\mathbb{N}$, $K(q^{2l})=K_l$.
\end{proposition}

\smallskip

 {\bf Proof.} The uniqueness of $K(u)$ is obvious. It follows from
(\ref{Kl}) that for $l=-1,-2,-3,\ldots$,
\begin{equation}\label{rec}K_l=\left(1+\sum \limits_{k=1}^m(-q^{2l})^k
\chi_k \right)\cdot K_{l+1}.
\end{equation}

 Remove the parentheses and reduce the right hand side of (\ref{Kl}) without
using any commutation relations (note that all the commutation relations
between the generators of $\mathbb{C}[\mathrm{Mat}_{mn}]_q^{\rm op}$ and
$\mathbb{C}[\overline{\mathrm{Mat}}_{mn}]_q$ are homogeneous of order two.
It suffices to prove that for each word over the alphabet
$\{\chi_1,\chi_2,\ldots,\chi_m \}$ the associated coefficient is a
polynomial of $u=q^{2l}$. To do this, observe that the coefficient at the
void word (free term) is 1, and the polynomial nature of other coefficients
is deducible via (\ref{rec}) using an induction argument with respect to the
length of the word. \hfill $\blacksquare$

\medskip

 Define the kernels $K_l \in \mathbb{C}[[\mathrm{Mat}_{mn}\times
\overline{\mathrm{Mat}}_{mn}]]_q$ for all $l \in \mathbb{C}$ by
$K_l=K(q^{2l})$, with $K(u)$ being the polynomial kernel whose existence and
uniqueness have just been proved.

\medskip

\begin{corollary} For $l \in \mathbb{N}$,
\begin{equation}\label{Kl1}K_l=\left(1+\sum \limits_{k=1}^m(-q^{2(l-1)})^k
\chi_k \right)^{-1}\cdot\left(1+\sum \limits_{k=1}^m(-q^{2(l-2)})^k \chi_k
\right)^{-1}\dots\left(1+\sum \limits_{k=1}^m(-1)^k \chi_k \right)^{-1}.
\end{equation}
\end{corollary}

\smallskip

 {\bf Proof.} The validity of this relation for $l \in \mathbb{C}$ follows
from the polynomial nature of $K(u)$ and the validity of (\ref{rec}) for $l
\in \{-1,-2, -3,\ldots \}$. \hfill $\blacksquare$

\medskip

 By a virtue of proposition \ref{fc}, for any $l \in \mathbb{C}$ one has a
well defined operator with kernel $K_l$:
$$\widehat{K}_l:D(\mathbb{U})_q \to D(\mathbb{U})_q';\qquad
\widehat{K}_l:{\rm id}\otimes \nu(K_l(1 \otimes f)y^l).$$

 Proposition \ref{plkl} admits the following generalization.

\medskip

\begin{proposition}\label{pilkl} For all $l \in \mathbb{C}$, $a \in U_q
\mathfrak{sl}_N$, one has the equality of operators from $D(\mathbb{U})_q$
to $D(\mathbb{U})_q'$:
$$\pi_l(a)\widehat{K}_l=\widehat{K}_l \pi_l(a).$$
\end{proposition}

\smallskip

 {\bf Proof.} It suffices to obtain the relation
\begin{equation}\label{wcr}\int
\limits_{\mathbb{U}_q}f_2(\pi_l(a)\widehat{K}_lf_1)d \nu=\int
\limits_{\mathbb{U}_q}f_2(\widehat{K}_l\pi_l(a)f_1)d \nu
\end{equation}
for all $f_1,f_2 \in D(\mathbb{U})_q$, $a \in U_q \mathfrak{sl}_N$, $l \in
\mathbb{C}$. By a virtue of proposition \ref{plkl} this relation is valid
for all $l \in-\mathbb{N}$. What remains is to prove that both hand sides of
\ref{wcr} are Laurent polynomials of the indeterminate $v=q^l$. As one can
observe from (\ref{ii}), (\ref{Thy}), it suffices to prove that the
integrands in (\ref{wcr}) are Laurent polynomials of $v=q^l$. (The function
$f(v)$ with values in $D(\mathbb{U})_q$ is called a Laurent polynomial if
all the operator valued functions $\Theta_{ij}(f(v))$, $i,j \in
\mathbb{Z}_+$, are Laurent polynomials (see section 4)). The polynomial
nature of the integrands in (\ref{wcr}) now follows from (\ref{pl1}) --
(\ref{pl3}) and proposition \ref{fc}. The latter proposition implies that
the formal series $K(u)=\sum \limits_{i,j=0}^\infty K^{(i,j)}$ which is
implicit in both hand sides of (\ref{wcr}) can be replaced by a finite sum
$\sum \limits_{i,j=0}^M K^{(i,j)}$, with $M=M(f_1,f_2)\in \mathbb{N}$.
Now what remains is to remind that the operator valued functions
$K^{(i,j)}(u)$, $u=q^{2l}$, are polynomials. \hfill $\blacksquare$

\bigskip

\section{q-analogues of Bergman kernels}

 In section 8 the kernels $K_\lambda \in \mathbb{C}[[\mathrm{Mat}_{mn}
\times \overline{\mathrm{Mat}}_{mn}]]_q$ have been defined in the special
case $\lambda \in \mathbb{Z}$ an explicit formula for $K_\lambda$ was
presented in section 8; the general case $\lambda \in \mathbb{C}$ is to be
considered in section 10.

 We are going to show that the orthogonal projections $P_ \lambda$ onto
Hardy-Bergman subspaces are integral operators with kernels $K_ \lambda$; in
different terms, these kernels are q-analogues of Bergman kernels.
(see \cite{FK}).

\medskip

\begin{theorem} For all $\lambda>N-1$, $f \in D(\mathbb{U})_q$ one has
\begin{equation}\label{il0}P_ \lambda f=({\rm id}\otimes \nu_ \lambda)(K_
\lambda(1 \otimes f)).
\end{equation}
\end{theorem}

\smallskip

 {\bf Proof.} Consider the special case $f=f_0$. Evidently,
\begin{equation}\label{il1}({\rm id}\otimes \nu_ \lambda)(K_ \lambda(1
\otimes f))=\left(\int \limits_{\mathbb{U}_q}f_0d \nu_ \lambda \right)\cdot
1.
\end{equation}
Prove that
\begin{equation}\label{il2}P_ \lambda
f_0=\left(\int\limits_{\mathbb{U}_q}f_0d \nu_ \lambda \right)\cdot 1.
\end{equation}
In fact, it follows from proposition \ref{Pl} that the element $P_ \lambda
f_0 \in \mathbb{C}[\mathrm{Mat}_{mn}]_q$ is subject to the relations $H_0(P_
\lambda f_0)=P_\lambda(H_0f_0)=0$. Hence $P_ \lambda f_0={\tt
const}(\lambda)\cdot 1$. On the other hand, $\|1 \|_ \lambda=1$ by a virtue
of proposition \ref{pnl}. What remains is to use the fact that $P_ \lambda$
is an orthogonal projection in $L^2(d \nu_ \lambda)_q$: ${\tt
const}(\lambda)=\dfrac{1}{\|1 \|_ \lambda}\displaystyle \int
\limits_{\mathbb{U}_q}1 \cdot f_0d \nu_ \lambda=\int
\limits_{\mathbb{U}_q}f_0d \nu_ \lambda$. Now (\ref{il1}), (\ref{il2}) imply
(\ref{il0}) for $f=f_0$.

 Our next step is to pass from the special case $f=f_0$ to the general case.
Let $L_ \lambda$ be the subspace of all those $f \in D(\mathbb{U})_q$ which
satisfy (\ref{il0}). Prove that $D(\mathbb{U})_q \subset L_ \lambda$.

 We know that $f_0 \in L_ \lambda$. By a virtue of propositions \ref{Pl} and
\ref{pilkl}, for all $a \in U_q \mathfrak{sl}_N$ one has $\pi_ \lambda(a)L_
\lambda \subset L_ \lambda$. Hence $\{\pi_ \lambda(a)f_0|\;a \in U_q
\mathfrak{sl}_N \}\subset L_ \lambda$. Apply proposition \ref{f0gdu} to
complete the proof. \hfill $\blacksquare$

\medskip

 Note that the measure $d \mu=d \nu_N$ is a q-analogue of the Lebesgue
measure in the matrix ball, and the kernel
$$K_N=\left(1+\sum_{k=1}^m\left(-q^{2(N-1)}\right)^k
\chi_k\right)^{-1}\left(1+\sum_{k=1}^m\left(-q^{2(N-2)}\right)^k
\chi_k\right)^{-1}\cdots\left(1+\sum_{k=1}^m(-1)^k \chi_k\right)^{-1}$$
 is a q-analogue of the ordinary Bergman kernel.

\medskip

 {\sc Remark 9.2.} The notation
$\mathbf{z}=(z_a^\alpha)_{a=1,\ldots,n;\;\alpha=1,\ldots,m}$,
$\boldsymbol{\zeta}=(\zeta_a^\alpha)_{a=1,\ldots,n;\;\alpha=1,\ldots,m}$,
allow one to rewrite (\ref{il0}) in a more appropriate form as
$$P_ \lambda=\int \limits_{\mathbb{U}_q}K_
\lambda(\mathbf{z},\boldsymbol{\zeta})f(\boldsymbol{\zeta})d
\nu_\lambda(\boldsymbol{\zeta}).$$

\bigskip

\section{Pairwise commuting kernels}

 Our purpose is to prove the following statements.

\medskip

\begin{lemma}\label{c1} In the algebra $\mathbb{C}[SU_m]_q^{\rm op}\otimes
\mathbb{C}[SU_m]_q$ the elements
\begin{equation}\label{his}\sum_{\genfrac{}{}{0mm}{1}{J',J''\subset\{1,
\ldots,m
\}}{{\rm card}(J')={\rm card}(J'')=k}}{z^{\wedge k}}_{J''}^{J'}\otimes
\left({z^{\wedge k}}_{J''}^{J'}\right)^\star,\qquad k=1,2,\ldots,m-1
\end{equation}
are pairwise commuting.
\end{lemma}

\medskip

\begin{lemma}\label{c2} In the algebra $\mathbb{C}[\mathrm{Mat}_{mm}]_q^{\rm
op}\otimes \mathbb{C}[\overline{\mathrm{Mat}}_{mm}]_q$ the elements
$\chi_k$, $k=1,\ldots,m$, given by (\ref{hik}), are pairwise commuting.
\end{lemma}

\medskip

\begin{proposition}\label{c3} In the algebra $\mathbb{C}[[\mathrm{Mat}_{mn}
\times
\overline{\mathrm{Mat}}_{mn}]]_q$ the elements $\chi_k$, $k=1,\ldots,m$,
given by (\ref{hik}), are pairwise commuting, and
\begin{equation}\label{Klam}K_ \lambda=\prod_{j=0}^\infty
\left(1+\sum_{k=1}^m(-q^{2(\lambda+j)})^k \chi_k
\right)\left(\prod_{j=0}^\infty \left(1+\sum_{k=1}^m(-q^{2j})^k \chi_k
\right)\right)^{-1}
\end{equation}
for all $\lambda \in \mathbb{C}$.
\end{proposition}

\medskip \stepcounter{theorem}

{\sc Remark 10.4.} The relations (\ref{Kl}) and (\ref{Kl1}) are special
cases of (\ref{Klam}).

\medskip

{\bf Proof} of lemma \ref{c1}. Apply the embedding of Hopf $*$-algebras
$\mathbb{C}[SU_m]_q\hookrightarrow(U_q \mathfrak{su}_m)^*$ (see
\cite{CP}).It sends the elements $z_a^\alpha \in \mathbb{C}[SU_m]_q$ to
matrix elements of operators of the vector representation $\pi$ of the
$*$-algebra $U_q \mathfrak{su}_m$ in the orthonormal basis of weight
vectors. The minors ${z^{\wedge k}}_{J''}^{J'}$ for $k>1$ could also be
treated in terms of exterior powers $\pi^{\wedge k}$of the vector
representation (see \cite{J}). Hence one has (see \cite{CP, S}):
\begin{equation}\label{mi1}\triangle({z^{\wedge
k}}_{J''}^{J'})=\sum_{{\rm card}(J)=k}{z^{\wedge k}}_{J}^{J'}\otimes
{z^{\wedge k}}_{J''}^{J},
\end{equation}
\begin{equation}\label{mi2}\left({z^{\wedge
k}}_{J'}^{J''}\right)^\star=S\left({z^{\wedge k}}_{J''}^{J'}\right),
\end{equation}
with $J,J',J''\subset \{1,2,\ldots,m \}$, $S$ and $\triangle$ being antipode
and comultiplication of the Hopf algebra $\mathbb{C}[SU_m]_q$. By a virtue
of (\ref{mi1}), (\ref{mi2}), the elements (\ref{his}) can be rewritten as
$\left((S \otimes{\rm id})\triangle\left(\displaystyle \sum
\limits_{\begin{array}{cc}\scriptstyle J \subset \{1,2,\ldots,m \}\\
\scriptstyle{\rm card}(J)=k \end{array}}{z^{\wedge
k}}_{J}^{J}\right)\right)^{\star \otimes \star}$. Since the linear map
$$\mathbb{C}[SU_m]_q \to \mathbb{C}[SU_m]_q^{\rm op}\otimes
\mathbb{C}[SU_m]_q,\qquad f \mapsto(S \otimes{\rm id})\triangle(f)$$ is a
homomorphism of algebras, it suffices to prove the pairwise commutativity of
the elements $\displaystyle \sum \limits_{\begin{array}{cc}\scriptstyle J
\subset \{1,2,\ldots,m \}\\ \scriptstyle{\rm card}(J)=k
\end{array}}{z^{\wedge k}}_{J}^{J}\subset \mathbb{C}[SU_m]_q$. What
remains is to apply the equivalence of the representations $\pi^{\wedge
k_1}\otimes \pi^{\wedge k_2}$ and $\pi^{\wedge k_2}\otimes \pi^{\wedge k_1}$
for all $k_1,k_2 \in \mathbb{Z}_+$. \hfill $\blacksquare$

\medskip

{\bf Proof} of lemma \ref{c2} can be obtained via replacing the quantum
groups $SU_m$ with the quantum group $U_m$ in the proof of lemma \ref{c1}.
Specifically, consider the Hopf $*$-algebra $\mathbb{C}[H]$ with the
standard comultiplication $\triangle$ and involution $\star$:
$\triangle(H)=H \otimes 1+1 \otimes H$, $H^\star=H$. Our definition implies
$U_q \mathfrak{u}_m=U_q \mathfrak{su}_m \otimes \mathbb{C}[H]$. What remains
is to demonstrate an embedding $\mathbb{C}[\mathrm{Mat}_{mm}]_q
\hookrightarrow(U_q \mathfrak{u}_m)^*$.

For that, consider the algebra $\mathbb{C}[GL_m]_q$ being a localization of
$\mathbb{C}[\mathrm{Mat}_{mm}]_q$ with respect to the multiplicative system
$(\det_q \mathbf{z})^\mathbb{N}$, where, as in an ordinary setting,
$$\det \nolimits_q \mathbf{z}=\sum_{s \in
S_m}(-q)^{l(s)}z_{s(1)}^1z_{s(2)}^2\ldots z_{s(m)}^m.$$

The algebra $\mathbb{C}[GL_m]_q$ is equipped with a structure of Hopf
algebra in a standard way and is called an algebra of functions on the
quantum group $GL_m$. Equip this Hopf algebra with an involution:
$$(z_a^\alpha)^*=(-q)^{a-\alpha}(\det \nolimits_q\mathbf{z})^{-1}\det
\nolimits_q(\mathbf{z}_a^\alpha),$$ with $\mathbf{z}_a^\alpha$ being thew
matrix derivable from $\mathbf{z}$ by discarding the line $\alpha$ and the
column $a$. The resulting Hopf $*$-algebra
$\mathbb{C}[U_m]_q=(\mathbb{C}[GL_m]_q,\star)$ will be called an algebra of
regular functions on the quantum group $U_m$. Now we have an embedding of
algebras $\mathbb{C}[\mathrm{Mat}_{mm}]_q \hookrightarrow \mathbb{C}[U_m]_q$
and an embedding $\mathbb{C}[U_m]_q \hookrightarrow (U_q \mathfrak{u}_m)^*$.
That is,
$$\mathbb{C}[\mathrm{Mat}_{mm}]_q \hookrightarrow \mathbb{C}[U_m]_q
\hookrightarrow(U_q \mathfrak{u}_m)^*.$$

Consider the involutive algebra $F=\mathbb{C}[SU_m]_q
\otimes\mathbb{C}[u,u^{-1}]$, $u^\star \overset{\rm def}{=}u^{-1}$. We need
embeddings of algebras
$$i_1:\mathbb{C}[\mathrm{Mat}_{mm}]_q^{\rm op}\hookrightarrow
F^{\rm op},\;i_2:\mathbb{C}[\overline{\mathrm{Mat}}_{mm}]_q \hookrightarrow
F,\;i:\mathbb{C}[\mathrm{Mat}_{mm}]_q^{\rm op} \otimes
\mathbb{C}[\overline{\mathrm{Mat}}_{mm}]_q \hookrightarrow F^{\rm op}\otimes
F,$$ given by
$$i_1(z_a^\alpha)=z_a^\alpha
\otimes u,\qquad i_2((z_a^\alpha)^*)=(z_a^\alpha)^\star \otimes
u^\star,\qquad i(f_1 \otimes f_2)=i_1(f_1)\otimes i_2(f_2).$$ It follows
from the definitions that
$$i(\chi_k)=\sum_{\begin{array}{c}\scriptstyle J',J''\subset \{1,\ldots,m
\}\\ \scriptstyle{\rm card}(J')={\rm card}(J'')=k \end{array}}({z^{\wedge
k}}_{J''}^{J'}\otimes u^k)\otimes(({z^{\wedge k}}_{J''}^{J'})^\star \otimes
u^{-k}),\qquad k \ne m,$$
$$i(\chi_m)=(1 \otimes u^m)\otimes(1 \otimes u^{-m}).$$
Thus we deduce from lemma \ref{c1} that
$$i(\chi_{k_1})i(\chi_{k_2})=i(\chi_{k_2})i(\chi_{k_1}),\qquad
k_1,k_2=1,\ldots,m.\eqno \blacksquare$$

\medskip

{\bf Proof} of proposition \ref{c3}. Show that the first statement reduces
to the special case $m=n$ which was considered in lemma \ref{c2}. Consider
the homomorphisms of algebras
$$j_1:\mathbb{C}[\mathrm{Mat}_{mn}]_q^{\rm op}\to
\mathbb{C}[\mathrm{Mat}_{mm}]_q^{\rm op},\qquad
j_2:\mathbb{C}[\overline{\mathrm{Mat}}_{mn}]_q \to
\mathbb{C}[\overline{\mathrm{Mat}}_{mm}]_q,$$
$$j:\mathbb{C}[\mathrm{Mat}_{mn}]_q^{\rm op}\otimes
\mathbb{C}[\overline{\mathrm{Mat}}_{mn}]_q \to
\mathbb{C}[\mathrm{Mat}_{mm}]_q^{\rm op}\otimes
\mathbb{C}[\overline{\mathrm{Mat}}_{mm}]_q$$ given by
$$j_1(z_a^\alpha)=\begin{cases}z_{a-(n-m)}^\alpha,& a>n-m \\0,& a
\le n-m\end{cases};\qquad j_2((z_a^\alpha)^*)=\begin{cases}(z_{a-(n-m)}^
\alpha)^*,&
a>n-m \\0,& a \le n-m\end{cases};$$
$$j(f_1 \otimes f_2)=j_1(f_1)\otimes j_2(f_2).$$
It suffices to prove the injectivity of the restriction of $j$ onto the
subalgebra $F_0$ generated by $\chi_k$, $k=1,2,\ldots,m$.

Let $\psi \in F_0$, $j(\psi)=0$. Choose $\lambda>N-1$ and consider the
integral operator with kernel $\psi$:
$$\widehat{\psi}:\mathbb{C}[\mathrm{Mat}_{mn}]_q \to
\mathbb{C}[\mathrm{Mat}_{mn}]_q;\qquad \widehat{\psi}:f \mapsto{\rm id}
\otimes \nu_ \lambda(\psi(1 \otimes f)).$$

Note that instead of the relation $\psi=0$ we may prove $\widehat{\psi}=0$
since the scalar product $(f_1,f_2)_ \lambda=\int
\limits_{\mathbb{U}_q}f_2^*f_1d \nu_ \lambda$ in the vector space
$\mathbb{C}[\mathrm{Mat}_{mn}]_q$ is nondegenerate.

Remind that $U_q \mathfrak{s}(\mathfrak{u}_n \times \mathfrak{u}_m)\subset
U_q \mathfrak{su}_{nm}$ is a Hopf $*$-subalgebra generated by $K_n^{\pm 1}$,
$\{E_j,F_j,K_j^{\pm 1}\}_{j \ne n}$. It is easy to show (see \cite{V}) that
the $U_q \mathfrak{s}(\mathfrak{u}_n \times \mathfrak{u}_m)$-invariance of
$\chi_k$, $k=1,\ldots,m$, implies the $U_q \mathfrak{s}(\mathfrak{u}_n
\times \mathfrak{u}_m)$-invariance of $\psi \in F_0$. Furthermore, the $U_q
\mathfrak{s}(\mathfrak{u}_n \times \mathfrak{u}_m)$-invariance of $y$
implies the $U_q \mathfrak{s}(\mathfrak{u}_n \times
\mathfrak{u}_m)$-invariance of the integral $\nu_ \lambda$. Hence the linear
map $\widehat{\psi}$ is a morphism of $U_q \mathfrak{s}(\mathfrak{u}_n
\times \mathfrak{u}_m)$-modules.

If $\psi \in F_0$ and $j(\psi)=0$, then one readily deduces that
$\widehat{\psi}$ is zero on the subalgebra generated by $z_a^\alpha$,
$a>n-m$. (In fact, if $\psi \in F_0$, then $j(\psi)=0$ is equivalent to ${\rm
id}\otimes j_2(\psi)=0$. Hence it suffices to prove the relation $\int
\limits_{\mathbb{U}_q}f^*\varphi d \nu_ \lambda=0$ for any element $\varphi$
of the subalgebra generated by $z_a^\alpha$, $a>n-m$, and any element $f$
such that $j_2(f^*)=0$. One can assume without loss of generality that
$$f=(z_1^1)^{k_{11}}(z_1^2)^{k_{21}}\ldots(z_1^m)^{k_{m1}}(z_2^1)^{k_{12}}
(z_2^2)^{k_{22}}\ldots(z_2^m)^{k_{m2}}\ldots(z_n^m)^{k_{mn}},$$
$$\psi=(z_{n-m+1}^1)^{l_{1,n-m+1}}(z_{n-m+1}^2)^{l_{2,n-m+1}}\ldots
(z_{n-m+1}^m)^{l_{m,n-m+1}}\ldots(z_n^m)^{l_{m,n}},$$ with $k_{i'j'} \ne 0$
for some $1 \leq i'\leq m$, $1 \leq j'\leq n-m$. However, in this case the
assumption $\int \limits_{\mathbb{U}_q}f^*\varphi d \nu_ \lambda \ne 0$
leads to a contradiction since for all $j=1,2,\ldots,n-1$ one has $H_j(f^*
\varphi)=0$, $\left(\sum \limits_{i=1}^mk_{ij}-\sum
\limits_{i=1}^ml_{ij}\right)-\left(\sum \limits_{i=1}^mk_{i,j+1}-\sum
\limits_{i=1}^ml_{i,j+1}\right)=0$,\footnote{We assume that $l_{ij}=0$ for
$j \le n-m$.} and that $H_0(f^*\varphi)=0$, $\sum
\limits_{j=1}^{n-1}\left(\sum \limits_{i=1}^mk_{ij}-\sum
\limits_{i=1}^ml_{ij}\right)=0$.)

The morphism of $U_q \mathfrak{s}(\mathfrak{u}_n \times
\mathfrak{u}_m)$-modules $\widehat{\psi}:\mathbb{C}[\mathrm{Mat}_{mn}]_q \to
\mathbb{C}[\mathrm{Mat}_{mn}]_q$ sends to zero all the generators
$$f_{j_1j_2 \ldots j_m}=\prod_{k=1}^m
\left({z^{\wedge k}}_{\{n-k+1,n-k+2,\ldots,n \}}^{\{1,2,\ldots,k
\}}\right)^{j_k}$$ of the $U_q \mathfrak{s}(\mathfrak{u}_n \times
\mathfrak{u}_m)$-module $\mathbb{C}[\mathrm{Mat}_{mn}]_q$. Hence
$\widehat{\psi}=0$, and thus the pairwise commutativity of $\chi_k$,
$k=1,\ldots,m$, is proved.

What remains is to obtain the relation (\ref{Klam}). Just the same argument
as that used in the proof of proposition \ref{pk}, allows one to establish
that the kernel $\prod \limits_{j=0}^\infty \left(1+\sum
\limits_{k=1}^m(-uq^{2j})^k\chi_k \right)$ is polynomial (see section 8).
Hence
$$K(u)=\prod \limits_{j=0}^\infty \left(1+\sum
\limits_{k=1}^m(-uq^{2j})^k\chi_k \right)\cdot \left(\prod
\limits_{j=0}^\infty \left(1+\sum \limits_{k=1}^m(-q^{2j})^k\chi_k
\right)\right)^{-1},$$ since the kernels in both hand sides are polynomials
and coincide with kernels (\ref{Kl}) as $u=q^{2l}$, $l \in-\mathbb{N}$. To
conclude, use the definition of $K_ \lambda$: $K_ \lambda=K(q^{2 \lambda})$.
\hfill $\blacksquare$

\bigskip

\section*{Appendix. Boundedness of the quantum matrix ball}

 Consider the faithful $*$-representation $\Pi$ of ${\rm Pol}({\rm
Mat}_{mn})_q$ in the {\sl pre-Hilbert} space $\widetilde{\cal H}$, described
in \cite[appendix 2]{SSV}. We use here the norm of the $m \times n$ matrix
with entries in ${\rm End}\,\widetilde{\cal H}$ defined as a norm of the
associated linear map $\bigoplus \limits_{a=1}^n \widetilde{\cal H}\to
\bigoplus \limits_{\alpha=1}^m \widetilde{\cal H}$.

\medskip

{\bf Proposition A.1.} {\it Let $Z$ and $\Pi(Z)$ be the matrices $(z_{\alpha
a})_{\alpha=1,\dots,m,\;a=1,\dots,n}$ and $(\Pi(z_{\alpha
a}))_{\alpha=1,\dots,m,\;a=1,\dots,n}$, respectively. Then $\|\Pi(Z)\|\le
1$.}

\smallskip

 {\bf Proof.} Let $S$ be the antipode of the Hopf algebra ${\mathbb
C}[SL_N]_q$; its action on the generators is given by a well known formula
(see \cite{CP}):
$$S(t_{a \beta})=(-q)^{a-\beta}\det \nolimits_q(T_{\beta a}),\qquad
a,\beta=1,\dots,N.\eqno(A.1)$$
(The matrix $T_{\beta a}$ in (A.1) is derived from
$T=(t_{ij})_{i,j=1,\dots,N}$ by discarding line $\beta$ and column $a$.)
Hence
$$\sum_{a=1}^N(-q)^{a-\beta}t_{\alpha a}\det
\nolimits_q(T_{\beta a})=\delta_{\alpha \beta},\qquad
\alpha,\beta=1,\dots,N,$$
or, equivalently,
$$-\sum_{c=1}^nt_{\alpha c}t_{\beta
c}^*+\sum_{\gamma=1}^mt_{\alpha,n+\gamma}t_{\beta,n+\gamma}^*=\delta_{\alpha
\beta},\qquad \alpha,\beta=1,\dots,m,$$
with $*$ being the involution in ${\rm Pol}(\widetilde{X})_q$ (see
\cite{SSV}). After introducing a notation
$$T_{11}=(t_{\alpha a})_{\alpha=1,\dots,m,\;a=1,\dots,n};\qquad
T_{12}=(t_{\alpha,n+\beta})_{\alpha,\beta=1,\dots,m},$$
$$T_{11}^*=(t_{a \alpha}^*)_{\alpha=1,\dots,m,\;a=1,\dots,n};\qquad
T_{12}^*=(t_{n+\beta,\alpha}^*)_{\alpha,\beta=1,\dots,m},$$
we get
$$-T_{11}T_{11}^*+T_{12}T_{12}^*=I.\eqno(A.2)$$

 It follows from (A.2) and (\ref{IZ}) that ${\EuScript
I}(I-ZZ^*)=T_{12}^{-1}(T_{12}^{-1})^*$.

 Apply the representation $\widetilde{\Pi}$ (see \cite{SSV}) to both parts
of the above relation. By a virtue of $\Pi=\widetilde{\Pi}{\EuScript I}$ we
obtain
$\Pi(I-ZZ^*)=\widetilde{\Pi}(T_{12}^{-1})\widetilde{\Pi}(T_{12}^{-1})^*\ge
0$. Hence $\Pi(Z)\Pi(Z)^* \le I$, $\|\Pi(Z)\|=\|\Pi(Z^*)\|\le 1$. \hfill
$\blacksquare$

\medskip

 Now a passage from $\widetilde{\cal H}$ to its completion allows one to
obtain a representation of the $*$-algebra ${\rm Pol}({\rm Mat}_{mn})_q$ by
bounded operators in a Hilbert space: $\|\Pi(z_a^\alpha)\|\le 1$,
$a=1,\dots,n$, $\alpha=1,\dots,m$.

\end{document}